\documentclass[12pt]{article}
\usepackage{epsf,epsfig,amsfonts,amsbsy,amssymb}

\input{diagrams}                  
\diagramstyle[height=2em,width=1em,scriptlabels,midshaft,
PS=dvips% noPS% 
]

\newarrow{Ratlongto} {}{dash}{}{dash}{->}
\newarrow{Into} C---{->}
\newarrow{To} ----{->}

\newtheorem{theo}{Theorem}
\newtheorem{prop}{Proposition}
\newtheorem{lemm}{Lemma}
\newtheorem{cor}{Corollary}

\def \PP {\mathbb{P}}
\def \QQ {\mathbb{Q}}
\def \NN {\mathbb{N}}
\def \ZZ {\mathbb{Z}}
\def \RR {\mathbb{R}}
\def \Chi {{\cal X}}
\def \longto {\longrightarrow}
\def \Tau {{\cal T}}

\def\pre{{\bf Proof : }}

\def\carre{\hfill $\square$}
\def\rem {\vspace{0.3cm}\noindent{\bf Remark : }}

\begin{document}
\begin{center}
\large {\bf
The GIT-equivalence \\
for $G$-line bundles\\}
\vspace{1cm}
N. Ressayre
\end{center}

\vspace{1.5cm}
{\bf Introduction}

Let $G$ be a reductive linear algebraic group acting algebraically on a 
projective variety $X$, both defined over an algebraically closed
field $k$.  
Let $L$ be an ample $G$-linearized line bundle on $X$. 
Geometric Invariant Theory (GIT) associates to $L$ a ``quotient''
$Y(L)$ of $X$ by $G$:
\begin{eqnarray}
\label{defY}
Y(L)=
\mbox{Proj}\left ( 
\bigoplus_{n\geq 0}\Gamma(X,L^{\otimes  n})^G
\right ).  
\end{eqnarray}

There is a natural $G$-invariant rational map  
$\pi\,:\,X--\to Y(L)$. 
The set where $\pi$ is defined is: 
$$
X^{\rm ss}(L)=\left \{
x\in X \mbox{ : }
\exists n>0\mbox{ and }
\sigma\in\Gamma(X,L^{\otimes n})^G \mbox{ such that }\sigma(x)\neq 0
\right \}.
$$
Points of $X^{\rm ss}(L)$ are said to be {\it semistable}\/ for $L$.
The map $\pi$ is obtained by gluing together  categorical quotients
of open affine $G$-stable subsets of $X^{\rm ss}(L)$.
As a consequence, each fiber of $\pi$ contains a unique closed
$G$-orbit in $X^{\rm ss}(L)$, and $(Y(L),\pi)$ is a categorical
quotient. 
In particular, $Y(L)$ only depends on $X^{\rm ss}(L)$ and is denoted by 
$X^{\rm ss}(L)//G$.
We also define  the open set  of {\it stable}\/ points:
$$
X^{\rm s}(L)=\left \{
x\in X^{\rm ss}(L) \mbox{ such that } G_x 
\mbox{ is finite and } 
G\cdot x \mbox{ is closed in } X^{\rm ss}(L)
\right \},
$$
where $G_x$ is the stabilizer of $x$.
It turns out that  $\pi^{-1}\left (\pi(x)\right )=G\cdot x$, 
for all $x\in X^{\rm s}(L)$.
We refer to  \cite{Kr} or \cite{GIT} for the classical properties
of this quotient.

Observe that this GIT-quotient is  not canonical:
it depends on a choice of an ample $G$-linearized line bundle $L$ 
over $X$ (sometimes called a polarization of $X$). 
During the last ten years, the question of variation of quotient
$X^{\rm ss}(L)//G$ under change of the ample $G-$linearized line
bundle $L$  has been an active research subject.
M. Brion and C. Procesi (see \cite{BP}) for a torus action and
M. Thaddeus (see  \cite{Th}) for the action of a reductive group
have studied the question when the linearization varies 
but the line bundle does not.
Recently, I. Dolgachev and Y. Hu obtained in \cite{DH} some results
when $L$ runs over all  ample $G$-linearized line bundles.
They introduced the notion of GIT-equivalence: 
two ample $G$-linearized line bundles $L_1$ and $L_2$ are said to be
{\it GIT-equivalent} if and only if $X^{\rm ss}(L_1)=X^{\rm ss}(L_2)$.
Using $G$-homological equivalence, they introduced a convex cone in a
finite dimensional real vector space: the {\it $G$-ample cone}.
The rational points of this cone give all GIT-quotients. 
The first theorem of \cite{DH} shows that only a finite number
of GIT-quotients can be obtained when $L$ varies. 
A second problem is to understand the geometry of the variation of
quotient. For this, given  three ample $G$-linearized
line bundles $L_-,\,L_0,\,L_+$ such that 
$X^{\rm ss}(L_-)\subset X^{\rm ss}(L_0)\supset X^{\rm ss}(L_+)$
one studies the transformation:

\begin{diagram}
  X^{\rm ss}(L_-)//G&                     &\rRatlongto&               &X^{\rm
    ss}(L_+)//G\\
&\rdTo_{\phi_-}&&\ldTo_{\phi_+}\\
&&X^{\rm ss}(L_0)//G
\end{diagram}
where the morphisms $\phi_\pm$ are induced by inclusions.
In \cite{Wa}, C. Walter gives a local model in etale topology for such a
transformation.

The goal of this paper is to study  the geometry of the
GIT-equivalence classes and its links with the inclusions $X^{\rm
  ss}(L_1)\subset X^{\rm ss}(L_2)$.
More precisely, if $C^G(X)$ denotes the $G$-ample cone, we show the
following: 

\vspace{0.2cm} 
\noindent
{\bf Theorem}
{\it Let $G$ be a reductive algebraic  group acting algebraically on a normal
projective variety $X$. Then:
\begin{enumerate}
\item For all $l_0\in C^G(X)$, 
$$C(l_0)=\left \{
l\in C^G(X) {\rm ~such~that~} X^{\rm ss}(l_0)\subset X^{\rm ss}(l)
\right \} $$
is a closed convex rational polyhedral cone in $C^G(X)$.
\item The cones $C(l)$ form a fan covering $C^G(X)$ (the {\it GIT-fan}).
\item The GIT-equivalence classes are the relative interior of these cones. 
\end{enumerate}
} 

\vspace{0.25cm}
We also study  the relations between the geometry of the
GIT-equivalence classes and the fibers of the morphisms 
$X^{\rm ss}(L_1)//G\longto X^{\rm ss}(L_2)//G$
induced by inclusions 
$X^{\rm ss}(L_1)\subset X^{\rm ss}(L_2)$.
We only use Geometric Invariant Theory contained in 
\cite{GIT} and some classical results on instability.
In particular, our arguments are valid over any algebraically closed
field and for any normal variety $X$.

We begin with recalling the numerical criterion of stability due to Mumford.
Then, we give some properties of this criterion due to Kempf,
Kirwan, Mumford, Ness, Hesselink etc. 
The first section ends with the finiteness theorem of 
I. Dolgachev and Y. Hu.
Following \cite{DH} and \cite{Th}, in a second section we introduce 
the algebraic
equivalence for $G$-linearized line bundles and the $G$-ample cone.   
The third section is devoted to the study of stability in the 
closure of an orbit. 
In the fourth and fifth sections, we study the geometry of the GIT-classes. 
We summarize our results by introducing the GIT-fan in Theorem
\ref{thGITfan}.  
Finally, we fix our attention on the fibers
of morphisms $\phi$ induced by inclusions of sets of semistable
points: we describe these fibers in the case where $G$ is a torus or
$\mbox{SL}(2)$ and we give examples for $G=k^*\times \mbox{SL}(2)$.

\vspace{0.3cm}
{\bf Acknowledgments.} I would like to thank I. Dolgachev and Y. Hu for
their questions which directed me and for their interest in my work.
Thanks also to C. Walter for useful discussions.  
I am especially grateful to Michel Brion who initiated me into 
Geometric Invariant Theory and who gave me precious advice
during the preparation of this article.

\section{The numerical criterion}

In this section we collect the notions and the results of
Geometric Invariant Theory (see \cite{GIT}) which will be used
throughout this paper.
We work over an algebraically closed field $k$ of arbitrary
characteristic. 

\subsection{The functions $\mu$}

Let $G$ be a reductive linear algebraic group acting
algebraically on a normal irreducible projective variety $X$.
Like in \cite{GIT}, we denote by  $\mbox{Pic}^G(X)$ the group of
$G$-linearized line bundles on $X$.
Let $L\in \mbox{Pic}^G(X)$.
Let $x$ be a point in $X$ and $\lambda$ be a one-parameter subgroup of
$G$.
Since $X$ is complete,  $\lim_{t \to 0}\lambda(t)x$ exists; let $x_0$
denote this limit.
The image of $\lambda$ fixes $x_0$ and so the group $k^\times$ acts
via $\lambda$ on the fiber $L_{x_0}$.
This action defines a character of $k^\times$, that is, an element of
$\ZZ$ denoted by $\mu^L(x,\lambda)$. 
The numbers $\mu^L(x,\lambda)$ satisfy the following properties:
\begin{enumerate}
\item $\mu^L(g\cdot x,g\cdot \lambda\cdot g^{-1})=\mu^L(x,\lambda)$ 
for any $g\in G$;
\item for fixed $x$ and $\lambda$, the map $L\mapsto \mu^L(x,\lambda)$
  is an homomorphism from $\mbox{Pic}^G(X)$ to $\ZZ$.
\end{enumerate}

The numbers $\mu^L(x,\lambda)$ are used in \cite{GIT} to give a 
numerical criterion for stability with respect to an ample $G-$linearized
line bundle $L$:
$$
\begin{array}{l}
  x\in X^{\rm ss}(L)\iff \mu^L(x,\lambda)\leq 0 \mbox{ for all
    one-parameter subgroups $\lambda$},\\
x\in X^{\rm s}(L)\iff \mu^L(x,\lambda)< 0 \mbox{ for all non trivial 
$\lambda$}.
\end{array}
$$

\subsection{Definition of the functions $\mbox{M}^\bullet(x)$}

Let $T$ be a maximal torus of $G$.
Denote the set of  one-parameter
subgroups of $T$ (resp. $G$) by $\Chi_*(T)$ (resp.  $\Chi_*(G)$).
We  denote the real vector space $\Chi_*(T)\otimes \RR$ by
$\Chi_*(T)_\RR$. 
The Weyl group $W$ of $T$ acts linearly on $\Chi_*(T)_\RR$.
Since $W$ is finite,
there exists a $W$-invariant Euclidean norm 
$\Vert\cdot \Vert$ on $\Chi_*(T)_\RR$. 
On the other hand, if $\lambda\in \Chi_*(G)$ there exists $g\in G$ such
that  $g\cdot \lambda\cdot g^{-1}\in\Chi_*(T)$. 
Moreover, if two elements of $\Chi_*(T)$
are conjugated by an element of $G$, then they are by an element of the
normalizer of $T$ (see \cite{Hum}).
This allows us to define the norm of $\lambda$ by
$\Vert\lambda\Vert=\Vert g\cdot \lambda \cdot g^{-1}\Vert$. 

Let $L\in \mbox{Pic}^G(X)$. We can now introduce the following notations:
$$
\begin{array}{cc}
\overline{\mu}^L(x,\lambda)=\frac{\mu^L(x,\lambda)}{\Vert
  \lambda\Vert},\hspace{0.2cm}&
\displaystyle{ \mbox{M}^L(x)=
\sup_{\lambda\in\Chi_*(G)}\overline{\mu}^L(x,\lambda).}
\end{array}
$$

Actually, it is shown in \cite{DH} that $\mbox{M}^L(x)$ is finite.
The functions $\mbox{M}^\bullet(x):\mbox{Pic}^G(X)\longto \RR$ 
will play a central role in the rest of the paper. 
 
\subsection{$\mbox{M}^\bullet(x)$ for a torus action}
\label{Mtore}

In this subsection we assume that $G=T$ is a torus.
Let $L$ be an ample $T$-linearized line bundle on $X$.
Then, there exist a $T$-module $V$ and a positive integer $n$ 
such that $X\subset\PP (V)$ and  
$L^{\otimes n}={\cal O}(1)_{\left|{ _X}\right.}$.
Replacing $L$ by $L^{\otimes n}$, we assume  that $n=1$.
We set $V_{\chi}=\{v\in V \mbox{ such that } \forall t\in
T~~t.v=\chi(t)v\}$. 
Then, we have:
$$
V=\bigoplus_{\mbox{\scriptsize $\chi\in \Chi^*(T)$}} V_{\chi}.
$$

Let $x\in X$ and $v\in V$ such that $[v]=x$.
There exist unique vectors $v_\chi\in V_\chi$ such that 
$v=\sum_{\chi}v_\chi$.
If $\lambda\in\Chi_*(T)$ then there exists, for all $\chi\in\Chi^*(T)$,
an  integer $<\lambda,\chi>$ such that for all  $t\in T$, we have:
$$
\lambda(t)\cdot v=\sum_{\chi\in\Chi^*(T)}t^{<\lambda,\chi>}v_\chi.
$$
We identify $\Chi^*(T)$ and $\Chi_*(T)$ with $\ZZ^r$ in such a way
that  $<\cdot ,\cdot >$ is the standard inner  product in $\ZZ^r$.

Set 
$\mbox{st}(x)=\{\chi\in\Chi^*(T) \mbox{ such that } v_\chi\neq 0\}$.
We have:
$$
\mu^L(x,\lambda)=\min_{\chi\in {\rm st}(x)} <\lambda,\chi>.
$$
 
The following proposition due to L.Ness gives a pleasant
interpretation of the number $\mbox{M}^L(x)$:

\begin{prop}(see \cite{Ne1})
\label{PropMtore}
With the above notations, we have:
\begin{enumerate}
\item
The point $x$ is unstable if and only if $0$ does not belong to the
convex hull of ${\rm st}(x)$.
In this case, ${\rm M}^L(x)$ is the distance from $0$ to this convex
hull.
\item
If $x$ is semistable, the opposite of ${\rm M}^L(x)$ is the 
distance from $0$ to the  boundary of this convex set.
\item
There exists $\lambda\in\Chi_*(T)$ such that
$\overline{\mu}^L(x,\lambda)={\rm M}^L(x)$.
If moreover $\lambda$ is indivisible, we call it an 
{\it adapted one-parameter subgroup for $x$}.
\item
If $x$ is unstable, there exists a unique adapted 
one-parameter subgroup for $x$.
\end{enumerate}
\end{prop}

\subsection{Properties of $M^\bullet(x)$}
\label{Mfini}

\begin{lemm}(see \cite{Ne1})
\label{MtoreMG} 
Let  $L$ be an ample G-linearized line bundle and  $T$ be a maximal
torus of  $G$. We denote by $r_T : {\rm Pic}^G(X)\to {\rm Pic}^T(X)$ 
the partial forgetful map.

Then, for all $x\in X$, the set of the numbers ${\rm M}^{r_T(L)}(g\cdot x)$ for  
$g\in G$ is finite and
${\rm M}^L(x)=\max_{g\in G} {\rm M}^{r_T(L)}(g\cdot x)$.
\end{lemm}

An indivisible one-parameter subgroup $\lambda$ of $G$ is said to be
{\it   adapted for   $x$ and $L$} if and only if
$\overline{\mu}^L(x,\lambda)=\mbox{M}^L(x)$.  
Denote by $\Lambda^L(x)$ the set of adapted one-parameter subgroups
for $x$.
 
\begin{cor}
\begin{enumerate}
\item The numbers $\mbox{M}^L(x)$ are finite (even if $L$ is not ample,
see  Proposition 1.1.6 in \cite{DH}).
\item If $L$ is ample, $\Lambda^L(x)$ is not empty.
\end{enumerate}
\end{cor}

Now, we can reformulate the  numerical criterion for stability: if 
$L$ is ample, we have
$$
\begin{array}{cc}
X^{\rm ss}(L)=\{x\in X : \mbox{M}^L(x)\leq 0\},&
X^{\rm s}(L)=\{x\in X : \mbox{M}^L(x)< 0\}.
\end{array}
$$

To a one-parameter subgroup $\lambda$ of $G$ we associate the
parabolic subgroup (see \cite{GIT}):
$$
P(\lambda)=\left \{
g\in G \mbox{ such that } 
\lim_{t\to 0}\lambda(t).g.\lambda(t)^{-1} 
\mbox{  exists in } G \right \}.
$$
Then, for $g\in P(\lambda)$, we have
$\mu^L(x,\lambda)=\mu^L(x,g.\lambda.g^{-1})$. 
The following theorem due to G. Kempf is a generalization of  the last
assertion of Proposition~\ref{PropMtore}.

\begin{theo}(see \cite{Ke})
\label{thKempf}
Let $x$ be an unstable point for an ample $G$-linearized line
bundle $L$. Then:
\begin{enumerate}
\item All the $P(\lambda)$ for $\lambda\in\Lambda^L(x)$ are equal. We
  denote by $P^L(x)$ this subgroup.
\item Any two elements of $\Lambda^L(x)$ are conjugated by an element
  of $P^L(x)$.   
\end{enumerate}
\end{theo}

We will also use the following theorem of L. Ness.

\begin{theo}(see \cite{Ne2})
\label{thNess}
Let $x$ and $L$ be as in the above theorem.
Let $\lambda$ be an adapted one-parameter subgroup for $x$ and $L$.
We consider $y=\lim_{t\to 0}\lambda(t)\cdot x$.
Then:
\begin{enumerate}
\item $\lambda\in\Lambda^L(y)$,
\item ${\rm M}^L(x)={\rm M}^L(y)$.
\end{enumerate}
\end{theo}

\subsection{Stratification of $X$ induced from $L$}

If $d>0$ and $<\tau>$ is a conjugacy class of one-parameter subgroups of
$G$, we set:
$$
S^L_{d,<\tau>}=\left \{ 
x\in X \mbox{ such that } \mbox{M}^L(x)=d \mbox{ and  } 
\Lambda^L(x)\cap<\tau>\neq\emptyset
\right \}.
$$
If ${\cal T}$ is the set of conjugacy classes of one-parameter subgroups,
the previous section gives us the following decomposition of $X$:
$$
X=X^{\rm ss}(L)\cup\bigcup_{d>0,\;<\tau>\in{\cal T}}S^L_{d,<\tau>}.
$$
W. Hesselink showed in \cite{He} that this union is a finite
stratification by $G$-stable locally closed subvarieties of $X$.
We will call it the {\it stratification induced from $L$}. 

Using this stratification, I. Dolgachev and Y. Hu have shown 
the following fundamental finiteness 
theorem (see Theorem 1.3.9 in \cite{DH}):

\begin{theo}
\label{thfini}
\begin{enumerate}
\item The set of locally closed subvarieties $S$ of $X$ which can be
  realized as the stratum $S^L_{d,<\tau>}$ for some ample
  $L\in{\rm Pic}^G(X),d>0$ and $\tau\in\Chi_*(G)$ is finite.
\item The set of possible open subsets of $X$ which can be realized as
  the set of semistable points with respect to some ample
  $G$-linearized line bundle is finite.
\end{enumerate}
\end{theo}
 
\section{The $G$-ample cone}
\subsection{Algebraic equivalence for $G$-line bundles}

Following I. Dolgachev and Y. Hu, we introduce the $G$-ample cone.
As M. Thaddeus in \cite{Th}, we use algebraic equivalence of $G$-line
bundles instead of homological equivalence.

Two elements $L_1$ and $L_2$ of $\mbox{Pic}^G(X)$ are said to be 
{\it $G$-algebraically equivalent} if there exist a connected
variety $S$, points $t_1,\,t_2\in S$, and a $G$-linearized line bundle
$L$ on $S\times X$ such that $L_{\left | \{t_1\}\times X\right .}=L_1$ and 
$L_{\left | \{t_2\}\times X\right .}=L_2$.
Here, $G$ acts on $S\times X$ via its action on the second factor.
Let $\mbox{NS}^G(X)$ denote the quotient of $\mbox{Pic}^G(X)$ by
this equivalence relation.
The following proposition is  analogous to Lemma 2.3.5 in \cite{DH}.

\begin{prop}
\label{propmuNS}
  Let $x\in X$ and $\lambda\in \Chi_*(G)$. Let
  $L_1,\,L_2\in{\rm Pic}^G(X)$ be $G$-algebraically equivalent.
Then, $\mu^{L_1}(x,\lambda)=\mu^{L_2}(x,\lambda)$.
\end{prop}

\pre
Let $S$, $L$, $t_1$ and $t_2$ be as in the definition of 
the $G$-algebraic equivalence. 
Denote by $y$ the point $\lim_{t\to 0}\lambda(t)\cdot x$.
The group $k^\times$ acts via $\lambda$  on the fibers $L_{(t,y)}$
for all $t\in S$. 
Consider the  map $S\longto \ZZ$, 
$t\mapsto \mu^L\left ( (t,y), \lambda \right )$. Obviously, this map
is locally constant. By connectness of $S$ we obtain
$\mu^{L_1}(x,\lambda)= \mu^L\left ( (t_1,y), \lambda \right )
=  \mu^L\left ( (t_2,y), \lambda \right )=
\mu^{L_2}(x,\lambda)$.
\carre

\vspace{0.3cm}
We denote by $\mbox{NS}^G(X)_\RR$ (resp. $\mbox{NS}^G(X)_\QQ$) 
the vector space $\mbox{NS}^G(X)\otimes\RR$ 
(resp. $\mbox{NS}^G(X)\otimes\QQ$).
A point $l\in \mbox{NS}^G(X)_\RR$ is said to be {\it rational} if it
belongs to $\mbox{NS}^G(X)_\QQ$.  
The forgetful map $\mbox{Pic}^G(X)\longto\mbox{Pic}(X)$ descends to 
$f:\mbox{NS}^G(X)\longto \mbox{NS}(X)$.
We denote by $\Chi^*(G)_\RR$ (resp. $\Chi^*(G)_\QQ$)
the vector space $\Chi^*(G)\otimes \RR$ (resp. $\Chi^*(G)\otimes
\QQ$). M. Thaddeus has proved in \cite{Th} the following:

\begin{prop}
  The map $f$ induces an exact sequence:
$$
0\longto \Chi^*(G)_\RR \longto {\rm NS}^G(X)_\RR\longto {\rm
  NS}(X)_\RR \longto 0
$$
\end{prop}

In particular, $\mbox{NS}^G(X)_\RR$ is a finite dimensional real
vector space.
Proposition~\ref{propmuNS} allows us to define
$\mu^\bullet(x,\lambda):\mbox{NS}^G(X)\longto\ZZ$ and 
$\mbox{M}^\bullet (x):\mbox{NS}^G(X)\longto\RR$.
By linearity, we can define 
$\mu^\bullet(x,\lambda):\mbox{NS}^G(X)_\RR\longto\RR$.
The function 
$\mbox{M}^\bullet (x):\mbox{NS}^G(X)_\RR\longto\RR$ 
is defined by the formula
$\mbox{M}^l(x)=\sup_\lambda\overline\mu^l(x,\lambda)$.

\begin{lemm}
Let $x\in X$. The function 
${\rm NS}^G(X)_\RR\longto\RR,\;l\mapsto {\rm M}^l(x)$ is convex and
positively homogeneous. In particular, it is continuous.
\end{lemm}

\pre
The functions $\mbox{M}^\bullet(x)$ are the suprema of a family of
linear forms. 
\carre
\vspace{0.6cm}

If $l\in \mbox{NS}^G(X)_\RR$, we set 
$$
\begin{array}{l}
X^{\rm ss}(l)=\{ 
x\in X \mbox{ such that } 
\mbox{M}^l(x)\leq 0\},\\ 
X^{\rm s}(l)=\{x\in X \mbox{ such that } \mbox{M}^l(x)< 0\}.
\end{array}
$$

We denote by $\mbox{NS}^G(X)^+_\RR$ the convex cone generated by the
classes of ample $G$-linearized line bundles in $\mbox{NS}^G(X)_\RR$.
This cone is open in $\mbox{NS}^G(X)_\RR$. 
Indeed, it is the preimage by $f$
of the ample cone, a strictly convex open cone in $\mbox{NS}(X)_\RR$.
A point $l$ in $\mbox{NS}^G(X)_\RR$ is said to be {\it ample} if
and only if it belongs to $\mbox{NS}^G(X)^+_\RR$.

A point $l\in \mbox{NS}^G(X)^+_\RR$ is said to be {\it effective} 
if and only if $X^{\rm ss}(l)$ is not empty. 
The set of effective points of $\mbox{NS}^G(X)^+_\RR$ is denoted by
$C^G(X)$ and called the {\it G-ample cone}. 
It is shown in \cite{DH} that the $G$-ample cone is convex.  
It may happen that $C^G(X)$ does not generate 
$\mbox{NS}^G(X)_\RR$; for example if $G$ is a product $H\times k^*$
and $k^*$ acts trivially on $X$. 

Two points $l$ and $l'$ in $C^G(X)$ are said to be {\it
  GIT-equivalent} if and only if $X^{\rm ss}(l)=X^{\rm ss}(l')$.
If $l_0\in \mbox{NS}^G(X)^+_\RR$, the {\it GIT-class} of $l_0$ is the set 
of all $l\in \mbox{NS}^G(X)^+_\RR$ that are GIT-equivalent to $l_0$.
The purpose of this paper is to describe the map 
$l\in C^G(X)\mapsto X^{\rm ss}(l)$. 
More precisely, we will describe the geometry of the GIT-classes.

\vspace{0.3cm}
\rem
Let $l\in C^G(X)$. If $l\in \mbox{NS}^G(X)_\QQ$, there
exists a $G$-linearized line bundle $L$ on $X$ such that 
$X^{\rm ss}(l)=X^{\rm ss}(L)$ and  $X^{\rm s}(l)=X^{\rm s}(L)$.
So, we have an algebraic quotient $X^{\rm ss}(l)//G$ with
all the properties of Geometric Invariant Theory. 
For example, if $l$ is rational then the sets $X^{\rm ss}(l)$ and
$X^{\rm s}(l)$ are open and a point $x\in X^{\rm ss}(l)$ belongs to 
$X^{\rm s}(l)$ if and only if its stabilizer is finite and its orbit
is closed in $X^{\rm ss}(l)$. 
At this step of the paper, we do not know if these properties hold for
$X^{\rm ss}(l)$ and $X^{\rm s}(l)$ when $l$ is not rational in $C^G(X)$.
This will often lead us to assume that: $l$ is rational.
Actually, this assumption will turn out to be unnecessary because
Proposition~\ref{PropGITouv} will show that any $l$ in $C^G(X)$ 
is GIT-equivalent to a rational one.

\subsection{A first property of the map $l\mapsto X^{\rm ss}(l)$}

The following proposition is a result of local monotonicity of the maps
$l\mapsto X^{\rm ss}(l)$ and $l\mapsto X^{\rm s}(l)$. 

\begin{prop}
\label{Propmon}
Let $l_0$ be a point in ${\rm NS}^G(X)^+_\RR$.
There exists a neighborhood $V$ of $l_0$ such that 
$X^{\rm s}(l_0)\subset  X^{\rm s}(l)\subset
X^{\rm ss}(l)\subset X^{\rm ss}(l_0)$
for all rational $l\in V$.
\end{prop}

\pre
By Theorem \ref{thfini}, there exist finitely  many open subsets 
$X^{\rm s}_1,\ldots,X^{\rm s}_{n}$ of $X$ 
such that $X^{\rm s}(l)$ is one of them for all
$l\in \mbox{NS}^G(X)_{\QQ}^+$. 
 We order these sets such that:
$$
\begin{array}{rl}
{\rm (i) }  &
X^{\rm s}(l_0)\not\subset X^{\rm s}_i \mbox{ for } i=1,\ldots,p\\
{\rm (ii) } &
X^{\rm s}(l_0)\subset X^{\rm s}_i \mbox{ for } i=p+1,\ldots,n.
\end{array}
$$

Let us fix some points $x_1,\ldots,x_p$ in  $X^{\rm s}(l_0)$  
such that $x_i\not\in X^{\rm s}_i$ for all $i=1,\ldots,p$.
By continuity of the functions $\mbox{M}^\bullet(x)$, there exist 
neighborhoods $V_{x_1},\ldots,V_{x_p}$ of $l_0$ in
$\mbox{NS}^G(X)^+_\RR$ such that:
$$
\forall l\in V_{x_i}\cap \mbox{NS}^G(X)_\QQ~~~x_i\in X^{\rm s}(l).
$$
Let $V$ denote the intersection of the $V_{x_i}$.
For all  $l\in V\cap \mbox{NS}^G(X)_{\QQ}$, $X^{\rm s}(l)$ is different from
$X^{\rm s}_i$ for all  $i=1,\ldots,p$.
Therefore, $X^{\rm s}(l)$ is one of the sets 
$X^{\rm s}_{p+1},\ldots,X^{\rm s}_{n}$ and contains $X^{\rm s}(l_0)$.

Let $(B_n)_{n\in\NN}$ be a sequence of balls centered at  $l_0$ and 
contained in $V$ such that the sequence of radii converges to zero. 
Let $X^{\rm ss}_1,\ldots,X^{\rm ss}_k$ be the open subsets of $X$
such that $X^{\rm ss}(l)$ is one of them for all
$l\in \mbox{NS}^G(X)_{\QQ}^+$. 
We set: 
$$
I_n=\left \{ i\;|\;
1\leq i\leq k ,\, \exists l\in B_n\cap \mbox{NS}^G(X)_{\QQ} 
\mbox{ such that }
l\neq l_0\mbox{ and } X^{\rm ss}(l)=X^{\rm ss}_i
\right \}
$$

The sequence $(I_n)$ is  decreasing. 
Since the sets $I_n$ are finite,
the sequence $I_n$ is  stationary (from a rank $N$).

Let $l\in B_N\cap \mbox{NS}^G(X)_{\QQ}$ and $x\in X^{\rm ss}(l)$.
For all $n$, there exists $l_n\in B_n$ such that
$x\in X^{\rm ss}(l_n)=X^{\rm ss}(l)$.
Since the sequence $l_n$ converges to $l_0$ and the function 
$\mbox{M}^\bullet(x)$ is continuous,  we have $\mbox{M}^{l_0}(x)\leq 0$; 
that is  $x\in X^{\rm ss}(l_0)$.
And so $X^{\rm ss}(l)\subset X^{\rm ss}(l_0)$.

\carre

\section{The stability set of a point}

\noindent
{\bf Definition}
If $x\in X$, the set  $\left \{ 
l\in \mbox{NS}^G(X)_\RR^+ \mbox{ such that } 
x\in X^{\rm ss}(l) \right \}$ 
is denoted  by $\Omega (x)$.
Observe that $\Omega(x)=\{l\in \mbox{NS}^G(X)_\RR^+ \mbox{ such that }
\mbox{M}^l(x)\leq 0\}$. 
The set $\Omega(x)$ is called {\it the stability set of $x$}.  

\vspace{0.3cm}
The aim of this section is to describe the geometry of 
the stability sets.
In other words, we fix a point $x\in X$ and study the semistability of
this point when $l$ varies in $\mbox{NS}^G(X)_\RR^+$.
Let $\overline{G\cdot x}$ denote the closure of the orbit of $x$.
In fact, the study of $\Omega (x)$ will lead us to consider the
GIT-classes of the $G$-variety $\overline{G\cdot x}$.

\subsection{A lemma on $\mu$}

\begin{lemm}
\label{Lemmulim}
  Let $l$ be a rational point in the $G$-ample cone.
Let $x\in X^{ss}(l)$ and $\lambda\in \Chi_*(G)$. 
We denote by $z$ the point $\lim_{t\to 0}\lambda(t)x$.

If $z$ is unstable for $l$ then $\mu^l(x,\lambda)<0$.
\end{lemm}

\pre
There exist a $G$-module $V$ and a positive integer $n$ such that $X$
can be embedded into $\PP(V)$ and the algebraic equivalence class of the 
restriction of ${\cal O}(1)$ to $X$ is  $l^{\otimes n}$. 
Replacing $l$ by $l^{\otimes n}$, we may assume that $n=1$.

Let us assume that $\mu^l(x,\lambda)$ is non negative.
Since $x$ is semistable for $l$,
we have $\mu^l(x,\lambda)=0$.
Let $T$ be a maximal torus in $G$ such that 
$\mbox{Im}(\lambda)\subset T$,
and let $g\in G$. 
We consider $\lambda'=g\lambda g^{-1},\,x'=g\cdot x$ and $z'=g\cdot z$. 
Then, $\mu^l(x',\lambda')= 0$.
We use the notations of Section~\ref{Mtore}.
Since $z'=\lim_{t\to 0} \lambda'(t)x'$,
$\mbox{st}(z')$ is equal to the intersection 
of $\mbox{st}(x')$ and  $\{\chi\mbox{ such  that } 
<\chi,\lambda'>=\mu^l(x',\lambda')\}$.
So the set $\mbox{st}(z')$ is a face of $\mbox{st}(x')$ containing $0$.
In particular, $z'$ is semistable for $r_T(l)$.
Finally, $z$ is semistable for $l$. 
This proves the lemma.
\carre

\subsection{A first result of rationality}

\begin{prop}
\label{PropOmRat}
Let $x\in X$. 
The stability set $\Omega(x)$ of $x$ is a convex cone and is closed in
${\rm NS}^G(X)^+_\RR$.
Moreover, the span of  $\Omega (x)$ is a rational vector 
subspace of ${\rm NS}^G(X)_\RR$.
In particular, $\Omega(x)$ is the closure of its rational points.   
\end{prop}

\pre
The first assertion  is obvious because the function $\mbox{M}^\bullet(x)$
is convex and positively homogeneous.
The last assertion is a direct consequence of the first ones. 
Let us prove the second one.
Let $F\subset\mbox{NS}^G(X)_\RR$ be the minimal rational 
vector subspace such that $\Omega (x)$ is contained in $F$.
Suppose that $\Omega (x)$ does not span $F$.
  
Since $\Omega (x)$ is convex, its interior as a subset of $F$ is empty.
This implies  that $\mbox{M}^l(x)=0$, for all $l\in\Omega(x)$.
Let $l$ be a point in $\Omega (x)$.
There exists a sequence $(l_n)$ of points in $F$ not in $\Omega (x)$ 
which converges to $l$.
Since $F$ is rational, we may assume that all $l_n$ are rational points.

By Theorem~\ref{thfini}, by extracting a subsequence, 
we may assume that all $l_n$ induce the same stratification $s$. 
For all $n$, $x$ is unstable for $l_n$. 
So there exists a non open stratum $S$ of $s$ containing $x$. 
Let $\lambda_0\in\Lambda^{l_0}(x)$ and $y=\lim_{t\to 0}\lambda_0(t)\cdot x$.
By Theorem~\ref{thNess}, $y\in S$, and so for all $n$ we have 
$\mbox{M}^{l_n}(x)=\mbox{M}^{l_n}(y)$.
But now, the continuity of the functions $\mbox{M}^\bullet (x)$ and 
$\mbox{M}^\bullet(y)$ implies that $\mbox{M}^l(y)=\mbox{M}^l(x)=0$.

Since $\lambda_0$ fixes $y$, 
we have $\mu^l(y,-\lambda_0)=-\mu^l(y,\lambda_0)$.
So, $\mbox{M}^l(y)=0$ implies $\mu^l(x,\lambda_0)=\mu^l(y,\lambda_0)=0$.
But $\mu^{l_0}(x,\lambda_0)>0$, so $\mu^\bullet (x,\lambda_0)$ is not
zero on $F$. 
Therefore, the equation $\mu^\bullet(x,\lambda_0)=0$ defines an
hyperplane of $F$ containing $l$. 

Moreover, the functions $\mu^\bullet (x,\lambda)$ are rational on
$\mbox{NS}^G(X)_\QQ$. 
In particular, the set of the functions equal to $\mu^\bullet
(x,\lambda)$  for some one-parameter subgroup $\lambda$ is countable.
Hence, $\Omega (x)$ is contained in a countable union of hyperplanes of $F$.
But $\Omega (x)$ is convex, 
thus $\Omega (x)$ is contained in such an hyperplane.
Therefore, $\Omega (x)$ is contained in a rational hyperplane of $F$.
Since $F$ is minimal, this is a contradiction.
The second assertion of the proposition is proved.
\carre

\vspace{0.3cm}
Note that Proposition~\ref{PropOmRat} is proved in \cite{DH} in the
special case when the codimension of $\Omega(x)$ is equal to one.

\begin{cor}
\label{CorOmRat}
The  number of stability sets is finite.
\end{cor}

\pre
It is clear that the stability set of each point of $X$ 
is an union of GIT-classes. 
So, by Theorem~\ref{thfini}, there exists only a finite number of sets of
the form $\Omega(x)\cap \mbox{NS}^G(X)_\QQ$ for $x\in X$.
Since  Proposition~\ref{PropOmRat} shows that  $\Omega(x)$ is the closure of 
$\Omega(x)\cap \mbox{NS}^G(X)_\QQ$,
the corollary is proved.
\carre

\vspace{0.3cm}
We also mention the well-known (see \cite{DH}) 

\begin{cor}
  The $G$-ample cone $C^G(X)$ is closed in ${\rm NS}^G(X)^+_\RR$.
\end{cor}

\pre
Since $C^G(X)=\cup_{x\in X}\Omega(x)$, Proposition~\ref{PropOmRat} and
Corollary~\ref{CorOmRat} imply that $C^G(X)$ is closed in
$\mbox{NS}^G(X)^+_\RR$. 
\carre

\subsection{Geometry of $\Omega(x)$}

The following lemma is essential in the study of the geometry of the
stability sets.

\begin{lemm}
\label{LemOmy}
  Let $x\in X$ and $z\in\overline{G\cdot x}-G\cdot x$. We assume that there exists
 a rational point $l_0\in C^G(x)$ such that $G\cdot z$ is closed in $X^{ss}(l_0)$.
Then, there exists $\lambda\in\Chi_*(G)$ such that
$$
\begin{array}{rl}
  (i)&\lim_{t\to 0}\lambda(t)\cdot x\in G\cdot z\\
  (ii)&
  \Omega (x)\subset \{l\in {\rm NS}^G(X)_\RR{\rm~such~that~}
\mu^l(x,\lambda)\leq 0\}\\
  (iii)&
  \Omega (z)=\{l\in {\rm NS}^G(X)_\RR {\rm~such~that~}\mu^l(x,\lambda)=0\} 
\cap \Omega (x)
\end{array}
$$
\end{lemm}

\pre
The Hilbert-Mumford theorem (see \cite{GIT}) applied to
$X^{ss}(l_0)$ gives us a 
$\lambda\in\Chi_*(G)$ such that $\lim_{t\to 0}\lambda(t)\cdot x\in G\cdot z$.
Denote by $z'$ this limit.
Since the image of $\lambda$ fixes $z'$, 
if $z'\in X^{\rm ss}(l)$ then $\mu^l(z',\lambda)$ and 
$\mu^l(z',-\lambda)=-\mu^l(z',\lambda)$ are negative or zero.
So we have 
$\Omega (z)=\Omega (z')\subset 
\{l\mbox{ such that }\mu^l(z',\lambda)=\mu^l(x,\lambda)=0\}$.

Let $l$ be a rational point in the stability set of $z$.
Since $X^{\rm ss}(l)$ is open 
and $z\in\overline{G\cdot x}$, we have $l\in \Omega (x)$.
But now, Proposition~\ref{PropOmRat} implies that  
$\Omega (z)$ is contained in $\Omega (x)$.
The inclusion 
$\Omega (x)\subset \{l\in {\rm NS}^G(X)_\RR\mbox{ such that
  }\mu^l(x,\lambda)\leq 0\}$ 
is  obvious; it implies that 
$\Omega (z)\subset \{l\in {\rm NS}^G(X)_\RR\mbox{ such that
  }\mu^l(x,\lambda)\leq 0\}\cap\Omega(x)$.
We show the opposite inclusion.

Let $l$ be a rational point in $\Omega (x)$
such that $\mu^l(x,\lambda)=0$.
Then, Lemma~\ref{Lemmulim} shows that $z$ is semistable for $l$.
So $l$ belongs to  the stability set of $z$.
We conclude by density of rational points in 
$\{l\mbox{ such that }\mu^l(x,\lambda)=0\}\cap \Omega (x)$
(which is the intersection of $\Omega(x)$ and a rational hyperplane). 
\carre

\vspace{0.3cm}
\noindent
{\bf Definition}
A {\it polyhedral cone} in $\mbox{NS}^G(X)^+_\RR$ is a subset of 
$\mbox{NS}^G(X)^+_\RR$ defined by a finite number of linear
inequalities.
Such a cone is said to be {\it rational} if the inequalities can be
chosen to be rational.
Let $C$ be a polyhedral cone in $\mbox{NS}^G(X)^+_\RR$. 
If $f$ is a linear form on $\mbox{NS}^G(X)_\RR$ non negative on $C$,  
the set of the points $c$ in $C$ such that $f(c)=0$ is said to be a
${\it face}$ of $C$.

\begin{prop}
\label{PropOmGeo}
Let $x\in X$.
\begin{enumerate}
\item
There exists $y\in \overline{G\cdot x}$ such that:
$l$ belongs to the relative interior of $\Omega(y)$ if and only if
$G\cdot y$ is closed in $X^{\rm ss}(l)$. 
Moreover, $\Omega(x)=\Omega(y)$. 
\item 
The stability set of $x$ is a convex rational polyhedral cone in 
${\rm NS}^G(X)^+_\RR$.
\item
The faces of $\Omega (x)$ are exactly the sets 
$\Omega (y)$ with $y\in\overline{G\cdot x}$.
\end{enumerate}
\end{prop}

\pre
Let $l_0$ be a rational point in the relative interior of $\Omega(x)$.
Let $y\in\overline{G\cdot x}$ such that $G\cdot y$ is closed in $X^{\rm ss}(l_0)$.
Lemma~\ref{LemOmy} shows that $\Omega(y)$ is a face of $\Omega(x)$. 
Since $l_0\in \Omega(y)$, we have $\Omega(x)=\Omega(y)$.
In particular, for all $z\in X$ there exists $z'\in\overline{G\cdot z}$
such that $\Omega(z)=\Omega(z')$ and $G\cdot z'$ is closed in 
$X^{\rm ss}(l)$ for some rational  point $l\in C^G(X)$.
Now, Lemma~\ref{LemOmy} shows that the sets $\Omega(z')$ with 
$z'\in\overline{G\cdot x}$ are faces of $\Omega(x)$.  

Let $l$ be a point in the relative boundary of $\Omega(y)$.
By Proposition~\ref{PropOmRat}, there exists  a sequence 
$(l_n)_{n\geq 1}$ 
of rational points in the vector space spanned by $\Omega(y)$, 
but out of $\Omega(y)$, which converges to $l$.
By Theorem~\ref{thfini}, by extracting a subsequence we may assume that
all $l_n$ induce the same stratification.
Now, like in the proof of Proposition~\ref{PropOmRat} we choose
$\lambda_l\in\Lambda^{l_1}(y)$ and we set
$z_l=\lim_{t\to 0}\lambda(t)\cdot y$. 
Then, $z_l\in\overline{G\cdot y}$ and $l\in\Omega(z_l)$.
Moreover, $\Omega(z_l)$ is contained in the hyperplane of
$\mbox{NS}^G(X)_\RR$ with equation $\mu^\bullet(y,\lambda_l)=0$, whereas
$\Omega(y)$ is not.
Thus $l$ belongs to $\Omega(z_l)$ which is a proper face of $\Omega(y)$.

We just proved that the relative boundary of $\Omega(y)$ is the union
of its faces $\Omega(z)$ for  some $z\in \overline{G\cdot y}$.
Now, Proposition~\ref{PropOmRat} and Corollary \ref{CorOmRat} 
imply the second assertion of the proposition.

Moreover, any face of codimension one of $\Omega(y)$ is equal to
$\Omega(z)$ for some $z\in \overline{G\cdot y}$.
By induction on the codimension of the face, this proves that  any
face of $\Omega(y)$ is equal to $\Omega(z)$ for some
$z\in\overline{G\cdot y}$.

Let us prove that $y$ satisfies the first assertion of the
proposition.
The above discussion shows that if $l$ belongs to the relative
boundary of $\Omega(y)$ then $G\cdot y$ is not closed in $X^{\rm ss}(l)$
(because $z_l\in(\overline{G\cdot y}-G\cdot y)\cap X^{\rm ss}(l)$).
Conversely, let $l'$ be a point in  $\Omega(y)$ such that $G\cdot y$ is
not closed in $X^{\rm ss}(l')$.
There exists $z'\in (\overline{G\cdot y}-G\cdot y)\cap X^{\rm ss}(l')$.
Since $G\cdot y$ is closed in $X^{\rm ss}(l_0)$, 
$z'\not\in X^{\rm ss}(l_0)$ and $\Omega(z')$ is a proper
face of $\Omega(y)$. 
Moreover,  $\Omega(z')$ contains $l'$.
Thus $l'$ does not belong to the relative interior of $\Omega(y)$. 
\carre

\vspace{0.3cm}
{\bf Definition}
A point $x$ is said to be {\it pivotal for $l$} if $x$ is semistable for $l$
(or for the GIT-class, $F$ of $l$) and $G\cdot x$ is closed in 
$X^{\rm  ss}(l)$. 
A point $x$ is said to be {\it pivotal for $\Omega$} if $\Omega$ is
the stability set of $x$ and $x$ is pivotal for some (or any) point in
the relative interior of $\Omega$.
  
\vspace{0.3cm}
Let us remark that I. Dolgachev and Y. Hu use in \cite{DH} a notion of
pivotal point which is close to  ours.

\section{Geometry of the GIT-classes}

\subsection{A fundamental lemma}

Up to now, we have studied the stability sets $\Omega(x)$ and for this the
quasi-homogeneous varieties $\overline{G\cdot x}$.
For an arbitrary $X$, the GIT-class of a point $l$ in 
$\mbox{NS}^G(X)^+_\RR$ depends on the relative positions of $l$ 
and the stability sets of the points of $X$.
To show how the various stability sets are related, the
main tool is the following:

\begin{lemm}
\label{LemXin}
Let $l_1$ and $l_2$ be two points in the
$G$-ample cone such that the set $X^{\rm ss}(l_1)$ is contained in 
$X^{\rm ss}(l_2)$.
Let $x$ be a pivotal point for $l_2$.

Then there exists a pivotal point $y$ for $l_1$ such that 
$x\in\overline{G\cdot y}$.
Moreover, $X^{\rm s}(l_2)$ is contained in $X^{\rm s}(l_1)$.
\end{lemm}

\pre
Let us consider the following commutative diagram:

\begin{diagram}
X^{\rm ss}(l_1)&\rInto^{\rm inclusion}&X^{\rm ss}(l_2) \\
\dTo_{\pi_{l_1}} &&\dTo_{\pi_{l_2}}\\
X^{\rm ss}(l_1)//G&\rTo^\phi&X^{\rm ss}(l_2)//G\\
\end{diagram}

The image of $\phi$ is equal to $\pi_{l_2}(X^{\rm ss}(l_1))$. 
Moreover, $\pi_{l_2}$ is surjective and $X^{\rm ss}(l_1)$ is dense in 
$X^{\rm ss}(l_2)$. So  $\phi$ is dominant. Since $X^{\rm ss}(l_1)//G$ is
complete, $\phi$ is surjective.

Therefore, there exists $y\in X^{\rm ss}(l_1)$ such that
$\phi(\pi_{l_1}(y))=\pi_{l_2}(x)$.
Since $y\in\pi_{l_2}^{-1}(\pi_{l_2}(x))$ and $G\cdot x$ is the unique closed
orbit in $\pi_{l_2}^{-1}(\pi_{l_2}(x))$, we have $x\in\overline{G\cdot y}$.
Moreover, $x\in\overline{G\cdot y'}$ for all 
$y'\in X^{\rm ss}(l_1)\cap\overline{G\cdot y}$.
Therefore we can find $y$ such that $G\cdot y$ is closed in 
$X^{\rm ss}(l_1)$ and $x\in\overline{G\cdot y}$.
This proves the first assertion of the lemma.

Let $x$ be a stable point for $l_2$ and let $y$ be as above. 
Since $G_x$ is finite, $G\cdot x=G\cdot y$ and $x$ is semistable for $l_1$.
Moreover, $G\cdot x$ is closed in $X^{\rm ss}(l_2)$ and hence in $X^{\rm
  ss}(l_1)$.
Thus, $x$ is stable for $l_1$.
\carre

\subsection{The geometry of a GIT-class}

Let $F$ be a GIT-class. 
We denote by $X^{\rm ss}(F)$ (resp. $X^{\rm s}(F)$) the subset
$X^{\rm ss}(l)$ (resp. $X^{\rm s}(l)$) of $X$ for some (or any) $l\in F$.  

\begin{lemm}
  \label{LemOmGIT}
Let $F$ be a GIT-class.
Then,
\[
F=\bigcap_{x\rm{~pivotal~for~}F} 
\rm{RelInt}(\Omega(x))
\]
where ${\rm RelInt}(\Omega(x))$ is the  relative interior of the
stability set of $x$. 
\end{lemm}

\pre
The first assertion of Proposition~\ref{PropOmGeo} shows that 
$F$ is contained in the intersection in the lemma. 
Conversely, let $l$ be a point in this intersection.
Then, $X^{\rm ss}(l)$ contains each closed orbit of $X^{\rm ss}(F)$,
so $X^{\rm ss}(F)$ is contained in $X^{\rm ss}(l)$.
Let $x\in X^{\rm ss}(l)$ such that $G\cdot x$ is closed in $X^{\rm ss}(l)$.
By Lemma~\ref{LemXin}, there exists $y\in X^{\rm ss}(F)$ such that 
$x\in \overline{G\cdot y}$ and $G\cdot y$ is closed in $X^{\rm ss}(F)$.
By assumption, $l$ belongs to the relative interior of the stability
set of $y$.
Then, Proposition~\ref{PropOmGeo} implies that $G\cdot y$ is closed in 
$X^{\rm  ss}(l)$. 
So $x\in G\cdot y$, and $X^{\rm ss}(l)=X^{\rm ss}(F)$.
The lemma is proved.
\carre

\begin{prop}
\label{PropGITouv}
Any GIT-class is the relative interior of a rational polyhedral 
cone in ${\rm NS}^G(X)_\RR^+$.
\end{prop}

\pre
The proposition follows immediately from Lemma~\ref{LemOmGIT}, 
Proposition~\ref{PropOmGeo} and Corollary~\ref{CorOmRat}.   
\carre

\vspace{0.3cm}
\noindent
{\bf Remark}
\begin{enumerate}
\item 
This proposition implies that any point in $C^G(X)$ is GIT-equivalent
to a rational point in $C^G(X)$. In particular, the second assertion of
Theorem~\ref{thfini}, Proposition~\ref{Propmon} and  Lemmas~\ref{Lemmulim}
and~\ref{LemXin} hold for any real point in $C^G(X)$.  
\item
  Let us assume for a moment that $X$ is a smooth complex variety.
Let $K$ be a maximal compact subgroup of $G$.
Let $\omega$ be a $K$-invariant K\"ahlerian symplectic form on $X$
and let $\phi$ be a moment map for the action of $K$.
In this situation, $X^{\rm ss}(\omega,\phi)=
\{x\in X\mbox{ such that } \overline{G\cdot x}\cap
\phi^{-1}(0)\neq\emptyset\}$
is called the set of semistable points. Then, $X^{\rm
  ss}(\omega,\phi)//G$ is a complex space, homeomorphic to
$\phi^{-1}(0)/K$. 
It turns out that $X^{\rm ss}(\omega,\phi)//G$ only depends on the
class of $\omega$ in $H^2(X,\RR)$ and on the choice of $\phi$ 
(see, for example Theorem 2.3.8 of \cite{DH}).
On the other hand, $\mbox{NS}(X)_\RR$ is a subspace of $H^2(X,\RR)$. 
The theory of Kempf-Ness shows that if $L$ and $\omega$ have the same
class in $H^2(X,\RR)$, then the complex spaces 
$X^{\rm ss}(\omega,\phi)//G$ and $X^{\rm ss}(L)//G$ are isomorphic.
So, the previous proposition  shows that if the class of $\omega$
belongs to $\mbox{NS}(X)_\RR$, the complex space $X^{\rm
  ss}(\omega,\phi)//G$ is a projective algebraic variety.
More generally, P. Heinzner and L. Migliorini showed in \cite{HeMi} that for all
$(\omega,\phi)$ the set $X^{\rm ss}(\omega,\phi)$ is equal to
$X^{\rm ss}(L)$ for some ample $G$-linearized line bundle $L$.
\end{enumerate}

\subsection{The inclusion relations $X^{\rm ss}(F)\subset X^{\rm ss}(F')$}

\begin{prop}
\label{PropGITin}
Let $F$ and $F'$ be two GIT-classes. The following assertions are
equivalent:
\begin{enumerate}
\item $F'$ intersects the closure of $F$ in ${\rm NS}^G(X)_\RR^+$
\item $F'$ is contained in the closure of $F$ in ${\rm NS}^G(X)_\RR^+$
\item $X^{\rm ss}(F)$ is contained in $X^{\rm ss}(F')$ 
\end{enumerate}
\end{prop}

\pre
It is sufficient to prove that 
$$
l\in \overline{F}\iff  X^{\rm ss}(F)\subset X^{\rm ss}(l),
$$
where $\overline{F}$ denotes the closure of $F$ in
$\mbox{NS}^G(X)^+_\RR$. 

If $l\in\overline{F}$, the continuity of the functions 
$\mbox{M}^\bullet(x)$ implies that $X^{\rm ss}(F)$ is contained in 
$X^{\rm ss}(l)$.
Conversely, let $l\in\mbox{NS}^G(X)^+_\RR$ be  such that 
$X^{\rm ss}(F)\subset X^{\rm ss}(l)$.  
Let $l_0$ be in $F$. 
Let $x$ be a pivotal point for $F$.
By Lemma~\ref{LemOmGIT}, $l_0$ belongs to the relative interior of
$\Omega(x)$. Moreover, the closure of $F$ is contained in $\Omega(x)$.
Then, the segment $[l_0;l[$ is contained in the relative interior of
$\Omega(x)$.
Now, Lemma~\ref{LemOmGIT} shows that $[l_0;l[$ is contained in $F$.
The proposition follows.
\carre

\begin{prop}
\label{PropGITfac}
The relative interior of a face of a GIT-class is a GIT-class.
\end{prop}

\pre
Let $F$ be a GIT-class.
By induction on the codimension of the face in $F$,
it is sufficient to prove the proposition for the maximal 
faces of $F$.
Let $F'$ be the relative interior of a maximal face of $F$.

Proposition~\ref{PropGITin} shows that the closure of $F$ is an 
union of GIT-classes. 
But, by Proposition~\ref{PropGITouv}, $F$ is the relative
interior of its closure.
Therefore, the relative boundary of $F$ is an union of GIT-classes. 
Since, by Proposition~\ref{PropGITouv} the GIT-classes are convex,
this implies that the closure of $F'$ is an union of GIT-classes.
Moreover, the GIT-classes are open in their closure.
So, $F'$ is an union of GIT-classes.
Thus, it is sufficient to prove that if $l_1$ and $l_2$ belong to
$F'$, then they are GIT-equivalent.

Let $y$ be a pivotal point for $l_1$.
Let us prove that $l_2$ belongs to $\Omega(y)$.
If  $y$ is semistable for $F$ then $l_2$ belongs to the closed cone
$\Omega(y)$.
Otherwise, by Lemma~\ref{LemXin} there exists a pivotal point $x$ for
$F$ such that $y\in\overline{G\cdot x}-G\cdot x$.
Indeed, by Proposition~\ref{PropGITin}, $X^{\rm ss}(F)$ is contained in 
$X^{\rm ss}(l_1)$. 
In particular, by Lemma~\ref{LemOmy}, there exists $\lambda\in\Chi_*(G)$
such that 
$\Omega(y)=\Omega(x)\cap \{l\in {\rm NS}^G(X)_\RR
{\rm~such~that~}\mu^l(x,\lambda)=0\}$ 
and, if $l$ belongs to
the relative interior of $\Omega(x)$ then $\mu^l(x,\lambda)<0$.
But, the first assertion of  Proposition~\ref{PropOmGeo}
shows that $F$ is contained in the relative interior of $\Omega(x)$. 
Thus, the intersection between the closure of $F$ and $\Omega(y)$  
is a proper face of $F$ containing $l_1$.
Then, $\overline{F}\cap \Omega (y)=\overline{F'}$.
In particular, $l_2$ belongs to $\Omega(y)$.

We just proved that all pivotal points for $l_1$ are semistable for
$l_2$.
This implies that $X^{\rm ss}(l_1)$ is contained in $X^{\rm ss}(l_2)$. 
By symmetry, $l_1$ and $l_2$ are GIT-equivalent.
The proposition is proved.
\carre

\section{Global geometry of the GIT-classes}

\subsection{A notion of wall and chamber}

In this section, we introduce a notion of wall and chamber.
Definitions close to ours have been considered by I. Dolgachev and
Y. Hu in \cite{DH} and by M. Thaddeus in \cite{Th}.

\vspace{0.3cm}
{\bf Definition}
A {\it wall} is a stability set of codimension one in $C^G(X)$.
A {\it chamber} is a GIT-class of codimension 0 in $C^G(X)$.
\vspace{0.4cm}

By continuity of the functions $\mbox{M}^\bullet(x)$, the GIT-classes
such that $X^{\rm ss}=X^{\rm s}$ are chambers. This is the notion of
chamber of \cite{DH}.
But all the chambers are not like that:
it is easy to find a $G$-action such that no GIT-class satisfies 
$X^{\rm ss}=X^{\rm s}$.
The appendix of \cite{DH} gives an example where
some chambers satisfy $X^{\rm ss}=X^{\rm s}$ and another does not. 

\begin{prop}
\label{LemMur}
The chambers are the connected components of the complement in
$C^G(X)$ of the union of the walls. The GIT-classes are the relative
interiors of the faces of the chambers.  
\end{prop}

\pre
The union of the closures of the chambers is closed. 
So its complement is an open subset of $C^G(X)$ covered
by the GIT-classes of codimension greater that one: this complement 
is empty. 
So each GIT-class meets the closure of a chamber. 
Proposition~\ref{PropGITfac} implies now the second assertion of
the proposition.
Moreover, each face of codimension one of a chamber is included in a
wall.
So, the $G$-ample cone is the disjoint union of the walls and the
chambers. The proposition follows immediately.   
\carre

\vspace{0.3cm}
Proposition~\ref{LemMur} implies that if we can determine the walls,  we
know all the GIT-classes. This remark is useful to calculate the
GIT-classes on examples.

\subsection{The GIT-fan}

\noindent
{\bf Definition} A {\it fan} $\Delta$ in $\mbox{NS}^G(X)^+_\RR$ is a
finite set of rational convex polyhedral cones in
$\mbox{NS}^G(X)^+_\RR$ such that
\begin{enumerate}
\item each face of a cone in $\Delta$ is also a cone in $\Delta$;
\item The intersection of two cones in $\Delta$ is a face of each of them.
\end{enumerate}

\vspace{0.3cm}
Most results of this paper about geometry of GIT-classes are
summarized in the following theorem announced in the introduction:

\begin{theo}
\label{thGITfan}
Let $G$ be a reductive algebraic group acting algebraically on a normal
projective variety $X$. Then:
\begin{enumerate}
\item For all $l_0\in C^G(X)$, 
$$C(l_0)=\left \{
l\in C^G(X) {\rm ~such~that~} X^{\rm ss}(l_0)\subset X^{\rm ss}(l)
\right \} $$
is a closed convex rational polyhedral cone in $C^G(X)$.
\item The cones $C(l)$ form a fan covering $C^G(X)$.
\item The GIT-classes are the relative interior of these cones. 
\end{enumerate}
This fan is called the {\it GIT-fan} for the action of $G$ on $X$.
\end{theo}

\pre
This theorem is a direct consequence of Propositions 
\ref{PropGITouv}, \ref{PropGITin}, \ref{PropGITfac} and~\ref{LemMur}.
\carre

\subsection{Existence of stable points}

The following proposition gives an easy criterion on $l$ for the
existence of stable points. It was first proved by I. Dolgachev and
Y.Hu with slightly different assumptions (see Propositions 3.2.8 and
3.3.5 in \cite{DH}).

\begin{prop}
\label{PropStabInt}
We assume that there exists $l_0\in C^G(X)$ such that $X^{\rm s}(l_0)$  is
not empty.
Then, for $l\in C^G(X)$, $X^{\rm s}(l)$ is not empty if and only
if $l$ belongs to the interior of $C^G(X)$.
\end{prop}

\pre
Let $l$ be a point in $C^G(X)$ such that $X^{\rm s}(l)$ is not empty.
Let $x\in X^{\rm s}(l)$.
By continuity, the function $M^\bullet(x)$ is negative on a
neighborhood of $l$.
Hence, $l$ belongs to the interior of $C^G(X)$.

Conversely, let $l$ be a point in the interior of $C^G(X)$.
There exists a point $l_1$ in $C^G(X)$ such that $l$ belongs to the
interval $]l_1;l_0]$.
Since the sets $X^{\rm ss}(l_1)$ and $X^{\rm s}(l_0)$ are open and non
empty, there exists $y\in X^{\rm ss}(l_1)\cap X^{\rm s}(l_0)$.
Then, by convexity of the function $M^\bullet(y)$, $y$ is stable for
$l$. 
The proposition is proved.
\carre

\section{The morphisms induced by inclusions}
\label{SectFibre}

Propositions~\ref{PropGITin} and~\ref{LemMur} show that for any
GIT-class $F$ there exists a chamber $C$ such that 
$X^{\rm ss}(C)\subset X^{\rm ss}(F)$. 
As a consequence, there exists a morphism 
$\phi\,:\, X^{\rm ss}(C)//G\longto X^{\rm ss}(F)//G$.
Moreover, in the proof of  Lemma~\ref{LemXin}, we already observed
that $\phi$ is surjective. 
So, the quotients corresponding to chambers
dominate the other quotients.

It could be interesting to compare two quotients corresponding to
two distinct chambers $C$ and $C'$. 
In this situation, there exists a sequence of chambers
$C=C_0,\,C_1,\,\ldots,\,C_m=C'$
such that for all $i=1,\ldots,m$, $C_{i-1}\cap C_i$ is a maximal
face of $C_{i-1}$ and $C_i$;
we denote by $F_i$ the relative interior of this face.
Then, by Theorem~\ref{thGITfan}, we have: 
$
X^{\rm ss}(C_0)\subset X^{\rm ss}(F_1)\supset 
X^{\rm ss}(C_1)\subset X^{\rm ss}(F_2)\supset\cdots 
\subset X^{\rm ss}(F_m)\supset X^{\rm ss}(C_m)
$.
These inclusions induce a sequence of morphisms:

\begin{equation}
\label{SuiteDeMorph}
\begin{diagram}
X^{\rm ss}(C_0)//G&                     &&                &X^{\rm ss}(C_1)//G
&                     &       &&\cdots   &  &  &X^{\rm ss}(C_m)//G.\\ 
                  &\rdTo(2,2)^{\phi_{-,1}}&      &\ldTo(2,2)^{\phi_{+,1}}&
&\rdTo(2,2)^{\phi_{-,2}}&       &         &  && \ldTo(2,2)^{\phi_{+,m}}\\
                  &                     &X^{\rm ss}(F_1)//G&        &
                      &&&\cdots &         &X^{\rm ss}(F_m)//G
\end{diagram}
\end{equation}

Observe that if there exists $l_0\in C^G(X)$ such that $X^{\rm
  s}(l_0)$ is not empty, then by Proposition~\ref{PropStabInt},
$X^{\rm s}(F_i)$ is not empty, for all $i=1,\ldots,m$.
In particular, the morphisms $\phi_{\pm,i}$ are birational.
These morphisms are studied in \cite{BP}, \cite{Th}, \cite{DH} and
\cite{Wa} in various degrees of generality.
In this section, our aim is to study the  geometry of
the fibers of the morphisms $\phi_{\pm,i}$.

More generally, let us fix  two effective ample
$G$-linearized line bundles $L_1$ and $L_2$ on $X$ such that 
$X^{\rm ss}(L_1)\subset X^{\rm ss}(L_2)$. 
Consider the following commutative diagram:

\begin{equation}
\label{diagram}  
\begin{diagram}
X^{\rm ss}(L_1) & \rInto^{\rm inclusion} & X^{\rm ss}(L_2) \\
\dTo_{\pi_{1}}& &\dTo_{\pi_{2}}\\ 
X^{\rm ss}(L_1)//G & \rTo^\phi& X^{\rm ss}(L_2)//G
\end{diagram}
\end{equation}

We will obtain general descriptions  of the fibers of $\phi$ 
in terms of quotients of affine subvarieties of $X$ by 
stabilizers of pivotal points for $L_2$. 
Then, we will come back to the situation in $(\ref{SuiteDeMorph})$ and
produce  pathological examples.  

\subsection{The fibers of the morphisms induced by inclusions}
\label{SecYS}

If $H$ is an algebraic group, we denote by $H^\circ$ the neutral
component of $H$. 
If $V$ is a $H$-module, the set of all $v\in V$ such that
$0\in\overline{H.v}$ is called the {\it nilcone of $V$}. 
Let $H$ be a reductive group and $Y$ be a $H$-variety.
If $\chi$ denotes a character of $H$, we denote by $L_\chi$ the
trivial line bundle on $Y$ linearized by:
$
h\cdot (y,t)=(h\cdot y,\chi(h)t) \hspace{0.3 cm}\forall 
h\in H,\,y\in Y {\rm ~and~}t\in k
$.
If $Y$ is affine, we denote by $Y//H$ the affine quotient, namely 
$\mbox{Spec}\left ( k[Y]^H\right )$.

\begin{lemm}
\label{lemS}
Let $H$ be a reductive group and $\Sigma$ be an affine $H$-variety (non
necessarily irreducible) with a fixed point $x$. 
We assume that $\{x\}$ is the only closed orbit of $H$ in $\Sigma$. 
Then:
\begin{enumerate}
\item The map $\Chi^*(H)\longto {\rm Pic}^H(\Sigma)$, 
$\chi \longmapsto L_\chi$ is an isomorphism.
\item We have  $\Sigma^{\rm ss}(L_\chi)=\Sigma$ if and only if $\chi$ 
is of finite order.
\end{enumerate}
\end{lemm}

\pre
Let $L$ be a linearized line bundle on $\Sigma$. The action of $H$ 
on the fiber $L_x$ gives a character $\chi$ of $H$. 
Since $\Sigma$ is affine, there exists  a section $s$ of $L$ such that
$s(x)\neq 0$. 
The set of zeroes of $s$ is closed, $H$-stable and does not contain $x$
: it is empty. So $L$ is trivial as a line bundle. This implies the first
assertion.

Let $L_\chi\in\mbox{Pic}^H(\Sigma)$ such that 
$\Sigma^{\rm ss}(L_\chi)=\Sigma$. 
By definition, for some $n>0$, there exists 
$s\in \Gamma(\Sigma,L_\chi^{\otimes n})^H$ 
such that $s(x)\neq 0$. 
This shows that $n\chi$ is trivial.
The second assertion of the lemma follows.
\carre

\begin{prop}
\label{PropYS}
Let $L_1,\,L_2,\,\pi_1,\,\pi_2$ and $\phi$ like in Diagram
\ref{diagram}.
Let $x$ be a pivotal point for $L_2$ which is not semistable for $L_1$. 
We denote by $H$ the  stabilizer of $x$.
Set $\Sigma=\{y\in X {\rm ~such~ that~ } x\in\overline{H\cdot y}\}$.
Then,
\begin{enumerate}
\item $H$ is reductive,
\item $\Sigma$ is affine,
\item $\phi^{-1}\left ( \pi_2(x)\right)$ is isomorphic to 
$\Sigma^{\rm ss}({L_1}_{\left | \Sigma \right .})//H$.
\end{enumerate}
Moreover, there exists a unique character $\chi$ of $H$ of infinite order 
such that ${L_1}_{\left | \Sigma \right .}=L_\chi$ in $\mbox{Pic}^H(\Sigma)$. 
If, in addition  $X$ is smooth then the $H$-variety $\Sigma$ is isomorphic
to the nilcone of a $H$-module.  
\end{prop}

\pre
Set $Z=\pi_2^{-1}\left (\pi_2(x)\right)$.
Then $Z$ is affine.
Since $G\cdot x$ is closed in $Z$, it is affine.
Hence, Matsushima's theorem shows that $H$ is reductive.

Let $\pi\,:\,Z\longto Z//H$ denote the quotient map.
Since $X^{\rm ss}(L_2)$ is open in $X$ and contains $x$, $\Sigma$ is
contained in $X^{\rm ss}(L_2)$. 
Then, $\Sigma$ is contained in $Z$ and 
$\pi^{-1}\left(\pi(x)\right )=\Sigma$.
In particular, $\Sigma$ is affine.  

Since $Z$ is affine and contains $G\cdot x$ as a unique closed orbit of
$G$, the Etale Slice Theorem of Luna (see \cite{GIT})
shows that $Z$ is equivariantly isomorphic to the fiber product
$G\times_H \Sigma$. 
But, by  the commutativity of Diagram~\ref{diagram} and the surjectivity 
of $\pi_1$, the fiber $\phi^{-1}\left(\pi_2(x)\right)$ is equal to
$
\pi_1\left (
X^{\rm ss}(L_1)\cap Z
\right )$.

Since $Z$ is $G$-stable and closed in $X^{\rm ss}(L_1)$, 
$
\pi_1\left (Z\cap X^{\rm ss}(L_1)\right )
$
is isomorphic to 
$
Z^{\rm ss}({L_1}_{\left | Z \right .})//G
$.
Moreover, for all $L\in\mbox{Pic}^G(G\times_H\Sigma)$ the restriction
map from $\Gamma(G\times_H\Sigma,L)^G$ to 
$\Gamma \left(
\Sigma,L_{\left |\Sigma \right .}
\right)^H$ is an isomorphism.
We conclude that:
$$
Y\simeq Z^{\rm ss}({L_1}_{\left | Z \right .})//G
\simeq \Sigma^{\rm ss}({L_1}_{\left | \Sigma \right .})//H.
$$

By Lemma \ref{lemS}, there exists a unique $\chi\in \Chi^*(H)$ such
that ${L_1}_{\left | \Sigma \right .}=L_\chi$.
Since $x\not\in \Sigma^{\rm ss}({L_1}_{\left | \Sigma \right .})$, 
Lemma \ref{lemS} shows that $\chi$ is of infinite order.

If in addition $X$ is smooth, the Etale Slice Theorem  shows that
$\Sigma$ is equivariantly isomorphic to the nilcone of the normal space
$T_xX/T_x(G\cdot x)$ to $G\cdot x$ in $X$ at the point $x$. 
The proposition is proved.
\carre

\vspace{0.3cm}
\rem
Conversely, let $H$ be a reductive group, $\chi$ be a character of
infinite order of $H$ and ${\cal N}$ be the nilcone of a $H$-module $V$.
Then, the variety ${\cal N}^{\rm ss}(L_\chi)//H$ occurs as a fiber of
a morphism $\phi$ like in Proposition \ref{PropYS}.

Indeed, let $X=\PP(V\oplus k)$. We define an action of $H$ on $X$ and
a linearization of ${\cal O}(1)$ by the formula:
$h\cdot (v,t)=(h\cdot v,t)$ for all $h\in H,\, v\in V$ and $t\in k$.
We denote by $L$ the line bundle ${\cal O}(1)$ so linearized.
Let $\pi\,:\,X^{\rm ss}(L)\longto X^{\rm ss}(L)//G$ be the quotient
map. 
It is easy to show that $\pi^{-1}\left (\pi([0:1])\right)={\cal  N}$.
On the other hand, by Proposition \ref{Propmon}, for $n$ large enough,
$X^{\rm ss}(L^{\otimes n}\otimes \chi)$ is contained in 
$X^{\rm ss}(L)$.
Let $\phi$ denote the induced morphism from 
$X^{\rm ss}(L^{\otimes n}\otimes \chi)//G$ to
$X^{\rm ss}(L)//G$.
Then, by Proposition \ref{PropYS}, 
$\phi^{-1}\left (\pi([0:1])\right )$ is isomorphic to 
${\cal N}^{\rm ss}(L_\chi)//H$. 

\begin{prop}
\label{PropYS2}
We assume that $G$ is connected.
We use the notations and assumptions of Proposition \ref{PropYS}.
Let $Y$ be an irreducible component of 
$\phi^{-1}\left(\pi_2(x)\right)$.
Then, there exists  an irreducible component $S$ of $\Sigma$ such
that $\pi_1\left( X^{\rm ss}(L_1)\cap \overline{G\cdot S}\right )$ is
equal to $Y$.

Let $H_S$ denote the stabilizer  of $S$ in $H$.
Let $\chi$ be as in Proposition \ref{PropYS} and 
$K_S= {\rm Ker}\chi\cap H_S$.
Then, we have:
  \begin{enumerate}
\item 
the group $H_S/{K_S}\simeq k^*$ acts on $S//{K_S}$
with a unique closed orbit which is a fixed point. 
\item
$ k[S//{K_S}]=\bigoplus_{n\in \NN}
k[S//{K_S}]_n$, where 
$$
k[S//{K_S}]_n=\left \{
f\in k[S] {\rm ~such~that~} h\cdot f=
\left( \chi(h)\right )^n f
\hspace{0.3cm}\forall h\in H_S
\right\}.
$$
\item 
$S^{\rm ss} ({L_1}_{\left | S \right .})//H_S$ is equal to 
$
{\rm Proj}\left ( 
\bigoplus_{n\geq 0} k[S//H_S]_n
\right ) 
$.  
  \end{enumerate}
Moreover, there exists a birational  finite morphism from 
$
S^{\rm ss} ({L_1}_{\left | S \right .})//H_S
$
onto $Y$. 
\end{prop}

\pre
Since by Proposition \ref{PropYS},
$\phi^{-1}\left(\pi_2(x)\right)\simeq 
\Sigma^{\rm ss}({L_1}_{\left|\Sigma\right.})//H$, 
the first assertion is obvious. 
Consider the surjective map:  
$$
\begin{array}{rccc}
\theta\,:&H\times_{H_S}S&\longto&H\cdot S\\
&(h:s)&\longmapsto&h\cdot s.
\end{array}
$$

Since $H^\circ$ stabilizes $S$, $H/H_S$ is finite.
In particular,  as a variety (without action of $H$), 
the fiber product $H\times_{H_S}S$ is isomorphic  to $H/H_S\times
S$. Thus, $\theta$ is finite.

Let $(H\cdot S)_0$ denote the set of all $x\in H\cdot S$ which belong to an
unique irreducible component of $H\cdot S$.
Set $S_0=(H\cdot S)_0\cap S$. 
One can easily prove that the restriction of $\theta$ to
$H\times_{H_S}S_0$ is an isomorphism onto $(H\cdot S)_0$. 
Thus, $\theta$ is finite, surjective and birational.

The restriction map $\rho \,:\, \mbox{Pic}^H(H\times_{H_S}S)
\longto \mbox{Pic}^{H_S}(S)$ is an isomorphism. 
We also denote  by ${L_1}_{\left | S \right .}$ the  $H$-linearized
line bundle $L$ on $H\times_{H_S}S$ such that  
$\rho(L)={L_1}_{\left | S \right .}$.
Then, $\theta$ induces a finite birational morphism:
$$
  \overline{\theta}\,:\,(H\times_{H_S}S)^{\rm ss} 
({L_1}_{\left | S \right .})//H \longto
(H\cdot S)^{\rm ss} 
({L_1}_{\left | H\cdot S\right .})//H\simeq Y.
$$

Since $G$ has a unique closed orbit in
$\pi_2^{-1}\left(\pi_2(x)\right)$,  the point $x$ is the
unique closed orbit of $H$ in $\Sigma$. 
As a consequence, $\{x\}$ is the unique closed orbit of $H_S$ in $S$.
Now, $k[S//{K_S}]^{H_S/{K_S}}=k[S]^{H_S}$ is equal to $k$.
As a consequence, $S//K_S$ contains a unique closed orbit
of $H_S/K_S$. 

Since $\chi$ is of infinite order,  the group $H_S/{K_S}$ is
isomorphic to $k^*$ via $\chi$. Then, the rational
$H_S/{K_S}$-module $ k[S//{K_S}]$ is equal to 
$\bigoplus_{n\in \ZZ} k[S//{K_S}]_n$.
Let $n$ be a negative integer and $f\in k[S//{K_S}]_n$.
By assumption, $S^{\rm ss}(L_\chi)=X^{\rm ss}(L_1)\cap S$ is not empty;
so there exists $n_0>0$ such that $k[S//{K_S}]_{n_0}$
contains a non zero function $f_0$.
Then $f^{n_0}f_0^n$ belongs to $k[S//{K_S}]_0=k[S]^{H_S}$ and must be
constant. By evaluating at $x$ this implies that $f=0$. 
Assertion $(ii)$ follows immediately.
 
Since $\Gamma (S,L_\chi^{\otimes n})=k[S//{K_S}]_n$ and 
${L_1}_{\left | S\right .}=L_\chi$, we have:
$$
(H\times_{H_S}S)^{\rm ss} 
({L_1}_{\left | S \right .})//H
   \simeq
S^{\rm ss} ({L_1}_{\left | S \right .})//H_S
   \simeq
{\rm Proj}\left ( 
\bigoplus_{n\geq 0} k[S//K_S]_n
\right ). 
$$
Now, the proposition follows from the properties of $\overline{\theta}$.
\carre

\vspace{0.3cm}
\rem
We will give an example where $\overline{\theta}$ is not an
isomorphism in Section \ref{SecExple}.

\subsection{The stabilizers of pivotal points and the stability sets}

\begin{prop}
\label{PropGlobOm}
Let $x\in X$. Then there exists a point $y\in X$ such that: 

\begin{enumerate}
\item 
$x\in\overline{G\cdot y}$
\item 
$\Omega(x)\subset \Omega(y)$
\item
 the interior of $\Omega(y)$ in $C^G(X)$ is not empty
\item
there exists a character $\chi$ of infinite order of $G_x$ such that
$G_y$ is contained in the kernel of $\chi$.  
\end{enumerate}
\end{prop}

\pre
Let us prove the three first assertion by induction on the codimension
of $\Omega(x)$. 
If this codimension is zero, then we can take $y=x$.
Otherwise,  let $l$ be a point in $\Omega(x)$.
Since the codimension of $\Omega(x)$ in $C^G(X)$ is at least one,
there exists a line ${\cal D}$ such that ${\cal D}\cap
\Omega(x)=\{l\}$ and ${\cal D}\cap C^G(X)\neq\{l\}$. 
Moreover, by Proposition~\ref{Propmon}, there exists $l'\in {\cal D}$ 
such that  $l'\neq l$ and $X^{\rm ss}(l')\subset X^{\rm ss}(l)$.
Lemma~\ref{LemXin} gives $x'\in X^{\rm ss}(l')$ such that 
$x\in \overline{G\cdot x'}$. 
But now, since $l'\not\in\Omega(x)$, $\Omega(x)$ is a proper face of 
$\Omega(x')$. By induction, the proposition holds for
$\Omega(x')$.
Therefore, there exists $y\in X$ which satisfies the three first
conditions of the proposition.

Let us prove that replacing $y$ by a point of $G\cdot y$, we can obtain the
last condition. 
Let $l_1$ (resp. $l_2$) be in the relative interior of $\Omega(y)$
(resp. $\Omega(x)$).
Consider the $G$-variety $\overline{G\cdot y}$ denoted by $X_y$. 
Since $G\cdot y$ is dense in $X_y$, the variety 
$
X_y^{\rm ss}({l_2}_{\left | X_y \right .})//G
$ is a point.
In particular, $X_y^{\rm ss}({l_2}_{\left | X_y \right .})$ is 
affine and contains $G\cdot x$ as unique closed orbit.
Let $\Sigma=\{z\in
X_y^{\rm ss}({l_2}_{\left | X_y \right .}) \mbox{ such that }
x\in\overline{G_x\cdot z} \}$.
By the Etale Slice Theorem, the $G$-variety 
$X_y^{\rm ss}({l_2}_{\left | X_y \right .})$ is isomorphic to
$G\times_{G_x}\Sigma$.
Replacing $y$ by a point in $G\cdot y$, we may assume that 
$y\in \Sigma$. 
Then, $G_y$ is contained in $G_x$.

Let $L_1$ be a $G$-linearized line bundle on $X$ in the homological 
class $l_1$.
Let $\chi$ be the character of $G_x$ which gives the action of $G_x$
on the fiber ${L_1}_x$.
Then, by Proposition \ref{PropYS}, the restriction of $L_1$ to
$\Sigma$ is $L_\chi$ and the order of $\chi$ is infinite. 
Moreover, $y$ is semistable for $L_1$. Thus, there exist
a positive integer $n$ and $f\in k[\Sigma]$ such that
$h\cdot f=\chi(h)^n f$ for all $h\in G_x$ and $f(y)\neq 0$.
In particular, the restriction of $\chi$ to $G_y$ is trivial.
The proposition is proved.
\carre

\begin{prop}
\label{PropRangCodim}
  Let $x$ be a point of $X$ which is pivotal for its stability set.
Then the rank of the character group of $G_x$ is at least the 
codimension of $\Omega(x)$ in the $G$-ample cone.
\end{prop}

\pre
Let $y$ be a point of $X$ which satisfies the conditions of
Proposition~\ref{PropGlobOm}.
Let $c$ denote the codimension of $\Omega(x)$ in $C^G(X)$.
There exists a sequence 
$\Omega(x)=\Omega_0\subset\Omega_1\subset\cdots\subset\Omega_c=\Omega(y)$
of faces of $\Omega(y)$ such that for all 
$i=1,\ldots,c$, the codimension of
$\Omega_{i-1}$ in $\Omega_{i}$ is equal to one.
For all $i=0,\ldots,c$, let $l_i$ be a point in the relative interior of 
$\Omega_i$.

Let $H$ denote the stabilizer of $x$ in $G$.
Consider $X_y=\overline{G\cdot y}$ and 
$\Sigma_y=\{z\in X_y \mbox{ such that } x\in \overline{H\cdot z}\}$.
By  Proposition~\ref{PropOmGeo}, for all $i=1,\ldots,c$, 
$X_y^{\rm ss}({l_i}_{\left | X_y\right.})$ is strictly contained in 
$X_y^{\rm ss}({l_{i-1}}_{\left | X_y\right.})$. 
Then, $X^{\rm ss}(l_0)\cap X_y\simeq G\times_H \Sigma_y$ implies that 
$\Sigma_y^{\rm ss}({l_i}_{\left | \Sigma_y\right.})$ 
is strictly contained in 
$\Sigma_y^{\rm ss}({l_{i-1}}_{\left | \Sigma_y\right.})$.

On the other hand, replacing $y$ by a point in $G\cdot y$, we may assume
the $y\in\Sigma_y$.
Let $S_y$ be an irreducible
component of $\Sigma_y$ containing $y$.
Since $X^{\rm ss}(l_0)\cap X_y$ is isomorphic to the fiber product 
$G\times_H \Sigma_y$, $G\times_H (H\cdot S_y)$ contains $G\cdot y$.
Then, $\Sigma_y=H\cdot S_y$ and thus, for all $i=0,\ldots,c$, we have 
$\Sigma_y^{\rm ss}({l_i}_{\left | \Sigma_y\right.})
=H\cdot 
S_y^{\rm ss}({l_i}_{\left | S_y\right.})
$.
In particular, 
$S_y^{\rm ss}({l_i}_{\left | S_y\right.})$ is strictly contained in 
$S_y^{\rm ss}({l_{i-1}}_{\left | S_y\right.})$,
for all $i=1,\ldots,c$. 

But by Theorem~\ref{thGITfan}, this implies that the 
$(l_i)_{0\leq  i\leq c}$ are affinely independent in 
$\mbox{NS}^{H^\circ}(S_y)_\RR$.
Moreover, by Proposition \ref{PropYS}, for all $i=0,\ldots,c$, there
exists a character $\chi_i$ of $G_x$ such that $L_{\chi_i}$ belongs to
the class ${l_i}_{\left | \Sigma_y\right.}$.
Thus, $\chi_0,\ldots,\chi_c$ are affinely independent.
The Proposition follows.
\carre

\rem 1. Proposition~\ref{PropRangCodim} is false with $\Omega(x)$
replaced by a GIT-class.
See for example: a linear  action of the two dimensional torus $T$ on
$\PP^3$  such that no  weight is  contained in the convex hull of the
three others.

2. If the codimension of $\Omega(x)$ is equal to the rank of $G$, 
Proposition~\ref{PropRangCodim} implies that the neutral component of 
the stabilizer of $x$ is a torus.
Arguing as in the proof of the last assertion of
Proposition~\ref{PropGlobOm}, this implies that for all $y\in X$ such
that $x\in\overline{G\cdot y}$, the neutral component of $G_y$ is a torus.

\subsection{The case where the stabilizer of a pivotal point is a torus}

\begin{prop}
\label{PropYStore}
  With the notations of Propositions~\ref{PropYS} and \ref{PropYS2},
if $H^\circ$ is a torus, then:
  \begin{enumerate}
  \item 
there exists  a one-parameter subgroup $\lambda$ of $H$ such
    that $S$ is an irreducible component of 
$\{y\in X : \lim_{t\to 0}\lambda(t)y=x\}$.
  \item
If in addition $X$ is smooth, then $S$ is equivariantly isomorphic to a
$H_S$-module and 
$
S^{\rm ss} ({L_1}_{\left | S \right .})//H_S
$
is the quotient of a projective toric variety by a finite group. 
Moreover, this quotient  is the normalization of $Y$.
  \item
If in addition the rank of $H^\circ$ is equal to one, then
$\phi^{-1}(\pi_2(x))$ is isomorphic to the quotient 
of a weighted projective space by a finite group.  
  \end{enumerate}
\end{prop}

\pre
The first assertion is a well-known fact about torus actions.
Let us assume that $X$ is smooth.
Then, $\Sigma$ is isomorphic to a nilcone for the action of the torus
$H^\circ$. 
In particular, the irreducible component $S$ of $\Sigma$ is 
equivariantly isomorphic to an $H_S$-module, denoted by $V$.

Let $\Tau$ be a maximal torus of $\mbox{GL}(V)$ containing the image of
$H^\circ$.
Let $L_V$ denote the restriction (after identification of $V$ with
$S$) of $L_1$ to $V$. 
The semistability of a point $v$ in $V$ for $L_V$ depends only on
the weights of $v$ for the action of $H^\circ$.
In particular, $\Tau$ stabilizes the open set $V^{\rm ss}(L_V)$.
Moreover, the actions of $H^\circ$ and of $\Tau$ on $V$ commute.
Hence, $\Tau$ acts on $V^{\rm ss}(L_V)//H^\circ$.
Because $\Tau$ has a dense orbit in $V$,
it has  a dense orbit in $V^{\rm ss}(L_V)//H^\circ$, too.
Thus, 
$
S^{\rm ss} ({L_1}_{\left | S \right .})//H_S
$ is the quotient of the projective toric variety 
$V^{\rm ss}(L_V)//H^\circ$ by the finite group $H_S/H^\circ$. 
Moreover, this quotient is normal. Then, the birational finite morphism of
Proposition \ref{PropYS2} is the normalization of $Y$.

From now on, we assume in addition that $H^\circ$ is
isomorphic to $k^*$.
Then, $S$ isomorphic to a $k^*$-module and by Proposition
\ref{PropYS2}, $
S^{\rm ss} ({L_1}_{\left | S \right .})//H_S
$ is the quotient of a weighted projective space by $H_S/H^\circ$. 
Moreover, an easy study of the linear actions of $k^*$ shows that 
$\Sigma$ has at most two irreducible components,
and that $\Sigma\cap V^{\rm ss}({L_1}_{\left | V\right . })$ is
contained in one irreducible component of $\Sigma$. 
In particular, with the notations of Proposition \ref{PropYS2},
$H_S=H$. 
But,  by Proposition \ref{PropYS}, $\phi^{-1}(\pi_2(x))$ is
isomorphic to $\Sigma^{\rm ss}({L_1}_{\left | \Sigma \right .})//H$.
Then, $\phi^{-1}(\pi_2(x))$ is isomorphic to 
$S^{\rm ss} ({L_1}_{\left | S \right .})//H$.
Assertion $(iii)$ of the proposition is proved.
\carre

\rem
With additional assumptions, Part $(iii)$ of Proposition
\ref{PropYStore} was proved, by different methods, by M. Thaddeus in
\cite{Th} and by I. Dolgachev and Y. Hu in \cite {DH}. 
In particular, I. Dolgachev and Y. Hu have shown that the assumptions of 
Proposition~\ref{PropYStore} are fulfilled for the diagonal action of
$G$ on $X\times G/B$, where $B$ is a Borel subgroup of $G$. 

\vspace{0.3cm}
The following proposition is an application of Propositions
\ref{PropRangCodim} and \ref{PropYStore} to the actions of $\mbox{SL}(2)$.

\begin{prop}
\label{PropWallSL2}
 Let $X$ be a ${\rm SL}(2)$-variety. 
Then:
  \begin{enumerate}
   \item 
Any wall is the intersection of an hyperplane and  
\vspace{0.3cm}${\rm  NS}^G(X)^+_\RR$. 
   \item
Let $\phi$ be a morphism like in Diagram~\ref{diagram}. 
We assume that $X$ is smooth. 
Then the fibers of $\phi$ are weighted projective spaces.  
  \end{enumerate}
\end{prop}

\pre
Since the rank of $\mbox{SL}(2)$ is equal to one, Proposition
\ref{PropRangCodim} shows that the codimension of $\Omega(x)$ is less
than one for all $x\in X$. 
Now, by Proposition~\ref{PropOmGeo}, a wall cannot have 
a boundary in ${\rm  NS}^G(X)^+_\RR$.
The first assertion follows immediately. 

From now on, we assume that $X$ is smooth.
Let $L_1,\,L_2,\,\pi_1,\,\pi_2,\,\phi,\,x,\,H$ and $\Sigma$ like in 
Propositions~\ref{PropYS}.  
Then, $H$ is a reductive subgroup of
$\mbox{SL}(2)$ which has a character of infinite order.
This implies easily that $H$ is a maximal torus of $\mbox{SL}(2)$.
Then, Propositions \ref{PropYS} and \ref{PropYS2} shows that 
$\phi^{-1}\left ( \pi_2(x)\right )$ is isomorphic to a weighted
projective space.
\carre

\vspace{0.3cm}
As an example, we refer to \cite{Po} for a detailed study of the
diagonal action of $\mbox{SL}(2)$ on $\left( \PP^1\right)^n$.

\subsection{Actions of $k^*\times \mbox{SL}(2)$}
\label{SecExple}

We use the notations of Proposition~\ref{PropYS}.
Let $F_1$ (resp. $F_2$) denote the GIT-class of $L_1$ (resp. $L_2$).
Proposition~\ref{PropYS} shows that the fibers of $\phi$ are  simpler
when $H$ is small.
Moreover, Propositions~\ref{PropGlobOm} and~\ref{PropRangCodim} show
that  $H$ is  small when  the dimension of $\Omega(x)$ is large.
As a consequence, natural restrictions in the study of the morphism
$\phi$ are:  

$$
\begin{array}{cl}
(H_1)& F_1 {\rm~is~a~chamber}\\
(H_2)&{\rm codim}(F_2)=1
\end{array}
$$ 
where $\mbox{codim}(F_2)$ is the codimension of $F_2$ in $C^G(X)$.
Moreover, these assumptions are fulfilled for the morphisms
$\phi_{\pm,i}$ considered in Diagram~\ref{SuiteDeMorph}.
On the other hand, Propositions~\ref{PropYS} and \ref{PropYStore} are
more precise if $X$ is smooth.

Actually, in \cite{DH} and \cite{Th} the assumption
$(H_1)$ is replaced by: $X^{\rm ss}(F_1)=X^{\rm s}(F_1)$.
Yet, if one wants to apply Construction \ref{SuiteDeMorph} to two
chambers $C$ and $C'$ such that  $X^{\rm  ss}(C)=X^{\rm s}(C)$ and 
$X^{\rm  ss}(C')=X^{\rm s}(C')$, one may have to consider
chambers $C_i$ such that $X^{\rm  ss}(C_i)\neq X^{\rm s}(C_i)$.
This happens indeed for $G=k^*\times \mbox{SL}(2)$, see the example in
the appendix of \cite{DH}.

On the other hand, if $G$ is a torus and $(H_1)$, $(H_2)$ hold, then
the fibers of $\phi$ are weighted projective spaces (this follows from
\cite{DH} or from Proposition \ref{PropYStore}).
The same holds if $G=\mbox{SL}(2)$ by Proposition~\ref{PropWallSL2}.
Yet, the following examples  show that for the actions of
$k^*\times\mbox{SL}(2)$ on a smooth variety, 
various varieties can occur as fibers of
$\phi$, even under the assumptions  $(H_1)$ and $(H_2)$.
We will give two examples where the morphism $\phi$ satisfies Assumptions
$(H_1)$ and $(H_2)$, whereas  $\phi$ has   
\begin{enumerate}
\item a reducible fiber, or
\item an irreducible  and non normal fiber.
\end{enumerate}

We will need the following technical lemma. 

\begin{lemm}
\label{lemtech}
Let $V$ be a finite dimensional vector space.
Let $V_+$ and $V_-$ be two vector subspaces. 
Let $H$ be a reductive group acting on $V_+\cup V_-$.
Let $H_+$ denote the stabilizer of $V_+$ in $H$.
We assume that $V_+\cup V_-=H\cdot V_+$ and that 
$H/H_{+}$ acts trivially on $(V_+\cap V_-)//H_+$.
Then, we have:
$$
(V_+\cup V_-)//H\simeq V_+//H_+.
$$ 
\end{lemm}

\pre
We claim that the restriction maps from $k[V_+\cup V_-]$ to $k[V_+]$
and $k[V_-]$ induces an isomorphism:
$$
k[V_+\cup V_-]\simeq \left \{
(f_+,f_-)\in k[V_+]\times k[V_-] \mbox{ such that }
{f_+}_{\left| V_+\cap V_- \right .}=
{f_-}_{\left| V_+\cap V_- \right .} \right \}.
$$
Endeed, let $f_\pm\in k[V_\pm]$ such that 
${f_+}_{\left| V_+\cap V_- \right .}=
{f_-}_{\left| V_+\cap V_- \right .} $.
Let $W$ be a vector subspace of $V$ and $W_\pm$ be two vector
subspaces of $V_\pm$ such that 
$V=W\oplus W_+\oplus W_-\oplus (V_+\cap V_-)$.
We define a function $\tilde{f}$ on $V$ by the formula:
 $\tilde{f}(w+w_++w_-+v)=f_+(w_++v)+f_-(w_-+v)-f_+(v)$ for all 
$w\in W,\, w_\pm\in W_\pm$ and $v\in V_+\cap V_-$.
Then, $\tilde{f}$ is regular on $V$ and the restrictions of
$\tilde{f}$ to $V_+$ and $V_-$ are respectively equal to $f_+$ and
$f_-$. The claim  follows easily.

Now, we consider the morphism 
$\theta : H\times_{H_+}V_+\longto V_+\cup V_-$ 
induced by the action of $H$ on $V_+\cup V_-$.
Via the comorphism $\theta^*$ of $\theta$, $k[V_+\cup V_-]$ is
identified to a subalgebra of $k[H\times_{H_+}V_+]$.
In particular, 
$\theta^*(k[V_+\cup V_-]^H)\subset 
k[H\times_{H_+}V_+]^H\simeq k[V_+]^{H_+}$. 

It is sufficient to prove that
$\theta^*\left(k[V_+\cup V_-]^H\right)=k[H\times_{H_+}V_+]^H$.
Let $f\in k[H\times_{H_+}V_+]^H$. 
Let $s\in V_+$ and $h\in H$ such that $h\cdot s\in V_+$.
Consider the quotient map, $\pi\,:\,V_+\longto V_+//H_+$.
If $h\in H_+$, then $\pi(s)=\pi(h\cdot s)$.
Otherwise, $h\cdot s\in V_+\cap V_-$.
But $H/H_{V_+}$ acts trivially on $(V_+\cap V_-)//H_+$.
Then, $\pi(s)=\pi(h\cdot s)$. 
Thus, for all $s\in V_+$ and $h\in H$ such that $h\cdot s\in V_+$, we have
$f(1:s)=f(1:h\cdot s)$.
Now, the claim implies that  $f$ belongs to 
$\theta^*\left(k[V_+\cup V_-]^H\right)$.
The lemma follows immediately.
\carre

\vspace{0.3cm}
\noindent
{\bf Examples }
We begin by fixing  some notation. 
From now on, $G$ denotes the group $k^*\times \mbox{SL}(2)$.
Let $\chi_0$ denote the character of $G$ defined by $\chi_0(t,g)=t$
for all $t\in k^*$ and $g\in \mbox{SL}(2)$.
Let $T$ be the maximal torus of $\mbox{SL}(2)$ of diagonal matrices and
$N(T)$ be its normalizer. 

Let $n\in\ZZ$ and $d\in \NN$. 
We denote by $V_d$ the $\mbox{SL}(2)$-module of binary forms of
degree $d$ in variables $a$ and $b$.  
We define an action of $k^*$ on $V_d$ which commutes to the action of
$\mbox{SL}(2)$ by the formula: 
$t\cdot v=t^nv$ for all $t\in k^*$ and $v\in V_d$. 
We obtain a $G$-module denoted by $V_{d,n}$.

\vspace{0.2cm}
Let $W$ be a $G$-module
(two specific choices of $W$ will be given below). Set  
$$
X=\PP(V_{2,0})\times \PP(W\oplus V_{0,1}).
$$
Let $\pi_1\,:\,X\longto \PP(V_{2,0})$ and $\pi_2\,:\,X\longto \PP(W\oplus
V_{0,1})$ denote the projection maps.
Set $L_1=\pi_1^*({\cal O}(1))$ and $L_2=\pi_2^*({\cal O}(1))$.
We linearize $L_1$  and $L_2$ canonically; that is, such that
$\Gamma(X,L_1)^*$  is $G$-isomorphic 
to  $V_{2,0}$ and $\Gamma(X,L_2)^*$ to $W\oplus V_{0,1}$. 
With the notations of Section \ref{SecYS}, we have : 
$\mbox{Pic}^G(X)=\mbox{NS}^G(X)=\ZZ L_1\oplus\ZZ L_2\oplus \ZZ L_{\chi_0}$.
Let $\RR_{>0}$ denote the interval $]0;+\infty[$.
Then, 
$$
\mbox{NS}^G(X)^+_\RR=\RR_{>0}L_1+\RR_{>0}L_2+\RR L_{\chi_0}.
$$

The decomposition of $V_2$ in eigenspaces for the action of
$T$ is $V_2=k.a^2\oplus k.ab\oplus k.b^2$. 
Denote the elements of $X$ by $([f],[(w,\tau)])$ where 
$f\in V_2,\,w\in W$ and $\tau\in V_{0,1}$.
Consider $$x_0=([ab],[(0,1)]).$$
We have 
$G_{x_0}=k^*\times N(T)$.
We identify, as vector spaces,  
$T_{x_0}X$ with $(k.a^2\oplus k.b^2)\times W$, 
and $T_{x_0}(G\cdot x_0)$ with $k.a^2\oplus k.b^2$.
In particular, $T_{x_0}X/T_{x_0}(G\cdot x_0)$ is isomorphic to $W$ as
a vector space.
Let $\rho\,:\,G\longto \mbox{GL}(W)$ denote the morphism induced by
the action of $G$ on $W$.  
We denote by  $W\otimes -\chi_0$ the representation of $k^*\times N(T)$
given by $\tilde{\rho}\,:\,k^*\times N(T)\longto \mbox{GL}(W)
,\; (t_0,t_1)\mapsto t_0^{-1}\rho((t_0,t_1))$. 
Then, we have the following isomorphism of $k^*\times N(T)$-modules:
$$
T_{x_0}X/T_{x_0}(G\cdot x_0)\simeq W\otimes -\chi_0.
$$

Let $m$ and $n$ be positive integers.
Let $V^{m,n}$ denote the $G$-module 
$
V_{2,0}^{\otimes m}\otimes (W\oplus V_{0,1})^{\otimes n}
$.
Consider the Segre embedding $i$ of $X$ in $\PP(V^{m,n})$.
Let $L_{m,n}$ be the line bundle ${\cal O}(1)$ on $\PP(V^{m,n})$
canonically linearized; in particular, 
$i^*(L_{m,n})=mL_1+nL_2$ in $\mbox{Pic}^G(X)$.

Let us use the notations of Section~\ref{Mtore} for the action of
$k^*\times T$ on $V^{m,n}$.
For all $x\in X$, the vertices of $\mbox{Conv(st(}i_{m,n}(x)))$ 
belong to the  set:
$$
{\cal V}=\{ m\alpha+n\beta \mbox{ such that } 
\alpha\in\mbox{st}(V_{2,0})\mbox{ and }
\beta\in\mbox{st}(W\oplus V_{0,1}) \}.
$$

Let $\chi_1$ denote the character of $k^*\times T$ defined by:
$\chi_1(t_0,t_1)=t_1$ for all $(t_0,t_1)\in k^*\times T$.
Then, $\Chi^*(k^*\times T)=\ZZ\chi_0\oplus\ZZ\chi_1$. 
We have: $\mbox{st(}i_{m,n}(x_0))=\{n\chi_0\}$. 
Moreover, for all $g\in G$, $n\chi_0$ belongs to
$\mbox{st}(i_{m,n}(g\cdot x_0))$. 
One can easily conclude that a point 
$L_1+q L_2+ \theta L_{\chi_0}\in \mbox{NS}^G(X)^+_\RR$ 
belongs to $\Omega(x_0)$ if and only if $q=\theta$.

\vspace{0.3cm}
{\bf A first choice of $W$} 

Set
$$
W=V_{1,-3}\oplus V_{1,-1}\oplus V_{1,1}\oplus V_{1,3}.
$$

\begin{figure}[htbp]
  \begin{center}
\leavevmode
\input{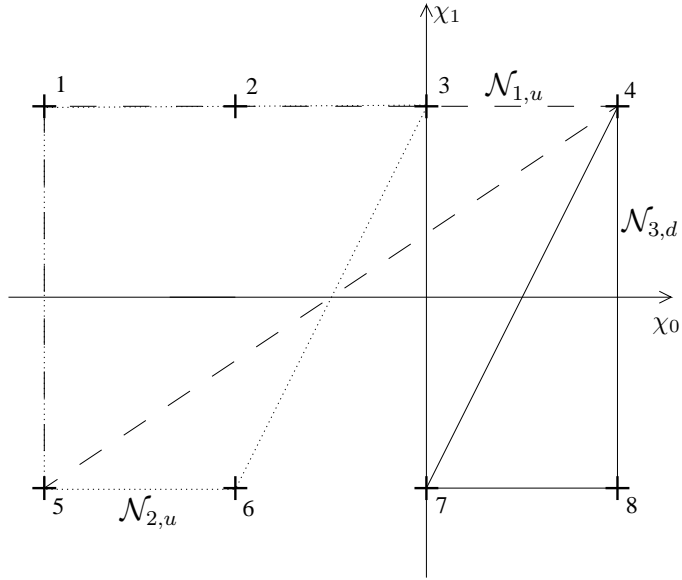}    
  \end{center}
\caption{Weights of $k^*\times T$ in $W\otimes -\chi_0$.}
\label{StateW}
%%% Magnification :  %%% 
\end{figure}

The weights of the action of $k^*\times T$ on $W\otimes -\chi_0$
are the crosses on Figure~\ref{StateW}. 
The meaning of the polytopes ${\cal N}_{1,u}$,
${\cal N}_{2,u}$ and ${\cal N}_{3,d}$ on Figure~\ref{StateW} 
will be explained later.

Set $q=n/m$.
The set ${\cal V}$ is represented by crosses on Figure~\ref{FigState}
after dilatation of ratio $1/m$. 
 
\begin{figure}[htbp]
  \begin{center}
\leavevmode
\input{state3.pstex_t}    
  \end{center}
\caption{The set ${\cal V}$.}
\label{FigState}
%%% Magnification : 50 %%% 
\end{figure}

Set $l_0=L_1+36(L_2+L_{\chi_0})$.
Note that $l_0$ belongs to $\Omega(x_0)$.
Set $l_-=l_0- L_{\chi_0}$ and $l_+=l_0+L_{\chi_0}$.
We denote respectively by $C_\pm$ the GIT-class of $l_\pm$.
One can easily see on Figure~\ref{FigState} that for $q=36$ the points
$35\chi_0$ and $37\chi_0$ do not belong to the boundary of a polytope
with vertices in ${\cal V}$.
We conclude that $C_+$ and $C_-$ are chambers such that 
$X^{\rm ss}(C_\pm)=X^{\rm s}(C_\pm)$.
By similar arguments one can show that the GIT-class $F_0$ of $l_0$ is
a face of $C_+$ and $C_-$.
Moreover, the codimension of $F_0$ in $C^G(X)$ is equal to one.
In particular, we have: 
$X^{\rm ss}(C_-)\subset X^{\rm ss}(F_0)\supset X^{\rm ss}(C_+)$.
These inclusions induce a diagram:

\begin{diagram}
X^{\rm ss}(C_-)//G&                     &  &                
&X^{\rm ss}(C_+)//G\\
                  &\rdTo(2,2)_{\phi_-}  &  &\ldTo(2,2)_{\phi_+}\\
                  &                     &X^{\rm ss}(F_0)//G
\end{diagram}

Let ${\cal N}$ denote the nilcone of
$W\otimes{-\chi_0}$ for the action of $k^*\times N(T)$.
Let $(\epsilon_1,\ldots,\epsilon_8)$ be a base of $W$ of eigenvectors 
for $k^*\times T$ such that the weight of $\epsilon_i$ is
the cross $i$ on Figure~\ref{StateW}.
Let $(x_1,\ldots,x_8)$ denote the dual basis of
$(\epsilon_1,\ldots,\epsilon_8)$. 
A vector $v\in W$ belongs to ${\cal N}$ if and only if $0$ does not
belong to the convex hull of the weights of $v$ for $k^*\times T$.
In particular, the irreducible components of ${\cal N}$ correspond to
the maximal convex hulls of weights of $k^*\times T$ which do not
contain $0$.
We obtain six irreducible components for ${\cal N}$. 
The equations of these six components are:
$$
\begin{array}{ll}
{\cal N}_{1,u}\,:\,x_6=x_7=x_8=0&
{\cal N}_{1,d}\,:\,x_2=x_3=x_4=0\\
{\cal N}_{2,u}\,:\,x_4=x_7=x_8=0&
{\cal N}_{2,d}\,:\,x_3=x_4=x_8=0\\
{\cal N}_{3,u}\,:\,x_1=x_2=x_5=x_6=x_7=0&
{\cal N}_{3,d}\,:\,x_1=x_2=x_3=x_5=x_6=0
\end{array}
$$

The convex hull of the weights of the action of $k^*\times T$ on 
${\cal N}_{1,u}$, ${\cal N}_{2,u}$ and ${\cal N}_{3,d}$ are
represented on Figure~\ref{StateW}.

Proposition~\ref{PropYS} shows that:    
\begin{eqnarray}
  \label{iso}
\phi_\pm^{-1}(\pi_0(x_0))\simeq
{\cal N}^{\rm ss}\left (L_{\pm \chi_0}\right )//(k^*\times N(T)).
\end{eqnarray}
In particular, 
$$
\phi_+^{-1}(\pi_0(x_0))
\simeq
({\cal N}_{3,u}\cup {\cal N}_{3,d})^{\rm ss}\left (L_{\chi_0}\right)
//(k^*\times N(T)).
$$
Moreover, by Lemma \ref{lemtech}, we have
$({\cal N}_{3,u}\cup {\cal N}_{3,d})//N(T)
\simeq {\cal N}_{3,u}//T$. 
This implies that 
$\phi_+^{-1}(\pi_0(x_0))\simeq 
{\cal N}^{\rm  ss}_{3,u}(L_{\chi_0})//(k^*\times T)$. 
In particular, this fiber is a projective toric variety of dimension
one. Then, we have  
$$
\phi_+^{-1}(\pi_0(x_0))\simeq \PP^1.
$$ 

By Isomorphism \ref{iso}, 
$({\cal N}_{1,u}\cup {\cal N}_{1,d})^{\rm ss}(L_{-\chi_0}) 
//(k^*\times N(T))$  is isomorphic to an irreducible component of 
$\phi_-^{-1}(\pi_0(x_0))$.
Write $k[{\cal N}_{1,u}]=k[x_1,\ldots,x_5]$.
The vector space $k[{\cal N}_{1,u}]^T$ is generated by the monomials
$x_1^{n_1}\cdots x_5^{n_5}$ such that $n_1+\cdots+n_4=n_5$.
Moreover, the weight of such a monomial for $k^*\times T$ is 
$-2(4n_1+3n+2n_3+n_4)\chi_0$.  
In particular,  the quotient
${\cal N}^{\rm ss}_{1,u}(L_{-\chi_0})//(k^*\times T)$
is isomorphic to the weighted projective space $\PP(1,2,3,4)$
and $N(T)/T$ acts trivially on 
$({\cal N}_{1,u}\cap {\cal N}_{1,d})//T$.
Then, Lemma \ref{lemtech} allows us to conclude that this irreducible
component of $\phi^{-1}_-\left(\pi_0(x_0)\right)$ is isomorphic to 
$\PP(1,2,3,4)$.

Let $\RR_{\geq 0}$ denote the set of  non negative real numbers.
If $d$ is a non negative integer, we denote by ${\cal P}_d$ the set of
$(n_1,n_2,n_3,n_5,n_6)\in \RR_{\geq 0}^5$ such that:
$$
\begin{array}[h]{rl}
&n_1+n_2+n_3=n_5+n_6\\
\mbox{and}&4(n_1+n_5)+2n_2+2n_6=d.
\end{array}
$$

Then, ${\cal P}_d$ is a convex polytope. 

One can  easily prove that 
${\cal N}^{\rm  ss}_{2,u}(L_{-\chi_0})//(k^*\times T)$ is 
isomorphic to  
$$
\mbox{Proj}\left (
\noindent\bigoplus_{d\geq 0}
\bigoplus_{(n_1,n_2,n_3,n_5,n_6) \in {\cal P}_d\cap\ZZ^5}
k.x_1^{n_1}x_2^{n_2}x_3^{n_3}x_5^{n_5}x_6^{n_6}
\right ).
$$  

Since the polytopes ${\cal P}_d$ are not simplicial, this toric
variety  is not a weighted projective space.  
Moreover, Isomorphism \ref{iso} and Lemma \ref{lemtech} imply that
an irreducible component of $\phi_-^{-1}(\pi_0(x_0))$
is isomorphic to this  toric variety. 

\vspace{0.3cm}
{\bf A second choice of $W$} 

Set
$$
W=V_{1,-1}\oplus V_{1,1}\oplus V_{3,3}.
$$

\begin{figure}[htbp]
  \begin{center}
\leavevmode
\input{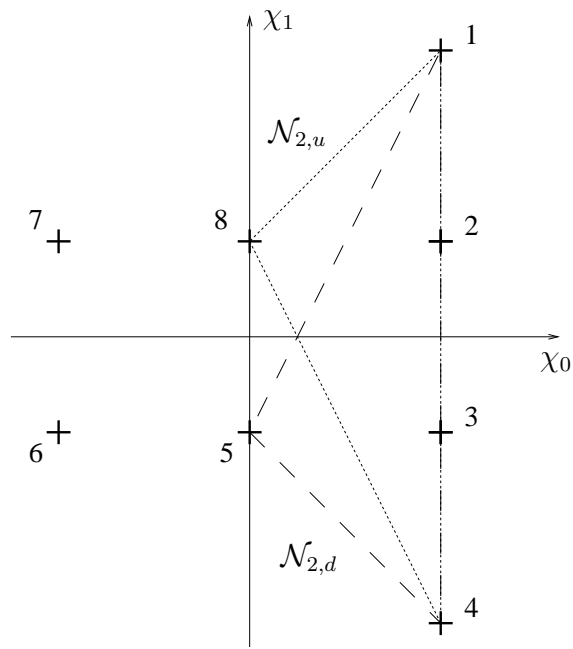}    
  \end{center}
\caption{Weights of $k^*\times T$ in $W\otimes -\chi_0$.}
\label{StateW2}
%%% Magnification : 50 %%% 
\end{figure}

The weights of the action of $k^*\times T$ on $W\otimes -\chi_0$
are the crosses on Figure~\ref{StateW2}. 

Like in the previous example, we set $l_0=L_1+36(L_2+L_{\chi_0})$,
$l_\pm=l_0\pm L_{\chi_0}$. 
The GIT-classes of $l_\pm$ are two chambers $C_\pm$ and the
GIT-class of $l_0$ is a maximal face $F_0$ of $C_+$ and $C_-$, too.
Let $\phi_+\,:\,X^{\rm ss}(C_+)//G\longto X^{\rm ss}(F_0)//G$ denote
the morphism induced by the inclusion 
$X^{\rm ss}(C_+)\subset X^{\rm ss}(F_0)$.

The nilcone ${\cal N}$ of $W\otimes-\chi_0$ has four irreducible
components. With obvious notations, the equations of these components
are:   
$$
\begin{array}{ll}
{\cal N}_{1,u}\,:\,x_1=x_2=x_3=x_4=x_5=0&
{\cal N}_{1,d}\,:\,x_1=x_2=x_3=x_4=x_8=0\\
{\cal N}_{2,u}\,:\,x_5=x_6=x_7=0&
{\cal N}_{2,d}\,:\,x_6=x_7=x_8=0
\end{array}
$$

By Proposition \ref{PropYS}, we have:
$$
\phi^{-1}_+\left(\pi_0(x_0)\right)
\simeq \left( {\cal N}_{2,u}\cup{\cal N}_{2,d}\right)^{\rm ss}
(L_{\chi_0})//\left(k^*\times N(T)\right). 
$$

Moreover, by Proposition \ref{PropYS2}, the natural map 
$$
\overline{\theta}\,:\,{\cal N}_{2,u}^{\rm ss}(L_{\chi_0})//(k^*\times T)
\longto \left( {\cal N}_{2,u}\cup{\cal N}_{2,d}\right)^{\rm ss}
(L_{\chi_0})//\left(k^*\times N(T)\right)
$$
is birational and finite.

The restriction of $\overline{\theta}$ to 
$\left( {\cal N}_{2,u}\cap{\cal N}_{2,d}\right)^{\rm ss}
(L_{\chi_0})//\left(k^*\times T\right)$ is the quotient by the action of
$N(T)/T$. 
But we have:
$$
k[{\cal N}_{2,u}\cap{\cal N}_{2,d}]^T
=k[x_1x_4,\,x_2x_3,\,x_1x_3^3,\,x_4x_2^3].
$$
Let $\alpha\,:\,k[X_1,\,X_2,\,X_3,\,X_4]\longto 
k[x_1x_4,\,x_2x_3,\,x_1x_3^3,\,x_4x_2^3]$ be the morphism defined by
$\alpha(X_1)=x_1x_4,\,\alpha(X_2)=x_2x_3,\,\alpha(X_3)=x_1x_3^3$ 
and $\alpha(X_4)=x_4x_2^3$.
Then, $\alpha$ induces an isomorphism $\bar{\alpha}$ from 
$k[X_1,\,X_2,\,X_3,\,X_4]/(X_3X_4-X_1X_2^3)$ onto
$k[{\cal N}_{2,u}\cap{\cal N}_{2,d}]^T$. 

Moreover, 
$$
k[{\cal N}_{2,u}\cap{\cal N}_{2,d}]^{N(T)}
= \bar{\alpha}\left( 
k[X_1,\,X_2,\,X_3+X_4]
\right ),
$$
which is strictly contained in $k[{\cal N}_{2,u}\cap{\cal N}_{2,d}]^T$.
In particular, $\overline{\theta}$ is not an isomorphism. 
Since, ${\cal N}_{2,u}^{\rm ss}(L_{\chi_0})//(k^*\times T)$ is normal, this
implies that $\phi^{-1}_+\left(\pi_0(x_0)\right)$ is irreducible but
not normal. 

\bibliographystyle{alpha}
\bibliography{biblio}

\begin{center}
  -\hspace{1em}$\diamondsuit$\hspace{1em}-
\end{center}

\noindent
Nicolas RESSAYRE\\
Universit\'e de Grenoble I\\
{\bf Institut Fourier}\\
UMR 5582 CNRS-UJF\\
UFR de Math\'ematiques\\
B.P 74\\
38402 SAINT MARTIN D'H\`ERES Cedex (France)\\
E-mail : {\tt ressayre@ujf.ujf-grenoble.fr}

\end{document}